\documentclass[12pt]{article}
\usepackage{amssymb}
\usepackage{amsbsy}
\usepackage{amsfonts}
\usepackage{amsmath}
\usepackage{amsopn}
\usepackage{latexsym}

\DeclareMathOperator{\e}{\mathbb{E}}

\DeclareMathOperator{\nr}{\mathbb{N}}

\DeclareMathOperator{\re}{\mathbb{R}}
\DeclareMathOperator{\p}{\mathbb{P}}

\newtheorem{theorem}{Theorem}

\newtheorem{corollary}[theorem]{Corollary}

\newtheorem{proposition}[theorem]{Proposition}

\begin{document}

\author{Jan Ob\l \'oj\thanks{\textit{e-mail: }\textrm{obloj@mimuw.edu.pl}}}\date{}\title{An explicit Skorokhod embedding for functionals of Markovian excursions}\maketitle 
\textit{Laboratoire de Probabilit\'es et Mod\`eles Al\'eatoires, Universit\'e Paris 6}\\\hspace*{2cm}\textit{ 4 pl. Jussieu - Bo\^{i}te 188, 75252 Paris Cedex 05, France.}

\textit{Faculty of Mathematics Warsaw University}\\\hspace*{2cm} \textit{Banacha 2, 02-097 Warszawa, Poland.}
\begin{abstract}We develop an explicit non-randomized solution to the Skorokhod embedding problem in an abstract setup of signed functionals of Markovian excursions. Our setting allows to solve the Skorokhod embedding problem, in particular, for diffusions and their (signed, scaled) age processes, for Az\'ema's martingale, for spectrally one-sided L\'evy processes and their reflected versions, for Bessel processes of dimension smaller than $2$, and for their age processes, as well as for the age process of excursions of Cox-Ingersoll-Ross processes.\\ This work is a continuation and an important generalization of Ob\l\'oj and Yor \cite{obloj_yor}. Our methodology, following \cite{obloj_yor}, is based on excursion theory and the solution to the Skorokhod embedding problem is described in terms of the It\^o measure of the functional. We also derive an embedding for positive functionals and we correct a mistake in the formula of Ob\l\'oj and Yor \cite{obloj_yor} for measures with atoms.\end{abstract} \noindent 2000 Mathematics Subject Classification: \textit{60G40, 60G44}\smallskip

\section{Introduction}
\label{sec:introch3}
The \emph{Skorokhod embedding problem} was first introduced and solved by Skorokhod \cite{MR32:3082b}, where it served to realize a random walk as a Brownian motion stopped at a sequence of stopping times. Since then, the problem has been generalized in a number of ways and has known many different solutions. At it simplest, the problem can be formulated as follows: given a Brownian motion $(B_t)$ and a centered probability measure $\mu$ with finite variance, find an integrable stopping time $T$ which embeds $\mu$ in $B$: $B_T\sim \mu$.
The condition of finite variance of $\mu$ was soon removed and the integrability of $T$ replaced with the requirement that $(B_{t\land T})$ is a uniformly integrable martingale. Unlike Skorokhod's original solution, subsequent solutions (e.g.\ Root \cite{MR38:6670}, Az\'ema and Yor \cite{MR82c:60073a}, Perkins \cite{MR88k:60085}, Jacka \cite{MR89j:60054}) were typically non-randomized and were often optimal in some sense. We refer to our survey paper \cite{genealogia} for further details.

Quite amazingly, this one problem continues to stimulate probabilists for over 40 years now and has actually seen a certain revival in the last years (e.g.\ Cox \cite{cox_thesis}, Cox and Hobson \cite{cox_hobson,cox_hobson_skew}, Ob\l\'oj \cite{genealogia,ja_maxprinciple}, Pistorious \cite{pistorius_skoro}). It also found new applications in the field of mathematical finance, such as pricing and hedging lookback and barrier options (Hobson \cite{hobson_lookback}, Brown, Hobson and Rogers \cite{MR1839367}). A new explicit solution, in discontinuous setup, was proposed in Ob\l\'oj and Yor \cite{obloj_yor}. Authors described an explicit and non-randomized solution to the Skorokhod embedding problem for the age of Brownian excursions, or more generally for positive functionals of Brownian excursions. However, they were only able to develop a randomized solution for the Az\'ema martingale.
Their work left two open challenges: firstly to extend the methodology to abstract Markovian setting, and secondly to extend the methodology to provide an explicit, non-randomized embedding for Az\'ema's martingale, or more generally for signed functionals of excursions.
The latter is very natural, as already argued by Ob\l\'oj and Yor \cite{obloj_yor}, and it actually motivated the present study. Indeed, Az\'ema's martingale -- the projection of Brownian motion on the filtration generated by the signs -- is an important process which, even though quite simple, inherits number of important properties from Brownian motion. It still finds various applications (e.g.~\c{C}etin et al.~\cite{MR2071419}) and so far no explicit non-randomized solution to the Skorokhod embedding problem for this process existed.

In this paper we solve both of the aforementioned open problems resulting from Ob\l\'oj and Yor \cite{obloj_yor}. We present an explicit non-randomized solution to the Skorokhod embedding problem for signed functionals of Markovian excursions. This abstract solution contains embeddings for Az\'ema's martingale, or signed age processes in general, for certain Bessel processes and for real-valued diffusions, to mention some examples. The stopping times we study here can be thought of as \emph{two-sided} generalizations of
the stopping times introduced by Ob\l\'oj and Yor \cite{obloj_yor} or, going back to the origins, of the stopping times introduced by Az\'ema and Yor \cite{MR82c:60073a}. We recall that Az\'ema and Yor \cite{MR82c:60073a} studied stopping times of the form
$T=\inf\{t\ge 0: B_t\ge \varphi(\sup_{s\le t}B_s)\}$, where $\varphi$ is an increasing function. In \cite{obloj_yor} Ob\l\'oj and Yor considered their analogue with the maximum process replaced with the local time at zero: $T=\inf\{t:F_t\ge \varphi(L_t)\}$, where $F_t$ was a positive functional of excursions (see Section \ref{sec:funcch3} below). The advantage of such stopping times is that $F_T=\varphi(L_T)$ and thus to describe the law of $F_T$ it suffices to describe the law of $L_T$.
Here we propose to investigate \emph{two-sided} stopping times, namely the stopping times of the form $T=\inf\{t: F_t\le -\varphi_-(L_t)\textrm{ or }F_t\ge \varphi_+(L_t)\}$.
Similar stopping times, in Brownian setting, were considered by Jeulin and Yor \cite{MR83c:60110} and Vallois \cite{MR86m:60200}, and also recently, in more generality, by Cox and Hobson \cite{cox_hobson_skew}. Given a probability measure $\mu$, we describe two increasing functions $\varphi_{+/-}$ such that $F_T$ has the distribution $\mu$. We consider functionals $F$ which are continuous over an excursion, except Section \ref{sec:disconch3} where discontinuous functionals are allowed to develop a solution to the embedding problem for spectrally one-sided L\'evy processes. As our stopping times allow to consider signed functionals $F$, a new difficulty arises, as compared with the study of Ob\l\'oj and Yor \cite{obloj_yor},
since $F_T$ can now take two values $-\varphi_-(L_T)$ or $\varphi_+(L_T)$. We will be able to deal with this new problem exploiting the properties of Poisson point processes. Similarly to \cite{obloj_yor}, theory of excursions and local times is the main tool in our work.

This paper is organized as follows. We first introduce the necessary notation and objects, in particular we discuss, in Section \ref{sec:excch3}, the excursion process of a Markov process and in Section \ref{sec:funcch3} we define the class of functionals we will consider and we clarify the terminology used throughout the paper. Then in Section \ref{sec:mainch3} we present our main results, in Theorem \ref{thm:mainch3} for signed functionals, and in Theorem \ref{thm:mainposch3} for positive functionals. The latter corrects a mistake found in the formula of Ob\l\'oj and Yor \cite{obloj_yor}. In subsequent two sections we develop applications of these results.
Section \ref{sec:azemach3} presents applications of Theorem \ref{thm:mainch3} and contains in particular an explicit solution to the Skorokhod embedding problem for the Az\'ema martingale, for Cauchy principal value associated with Brownian local times (over one excursion), for skew Brownian motion and for Brownian motion itself.
Section \ref{sec:pozch3} contains applications of Theorem \ref{thm:mainposch3} and develops explicit solutions to the Skorokhod embedding problem for Bessel processes of dimension $\delta\in (0,2)$ and their age processes, as well as for the age process of excursions of Cox-Ingersoll-Ross processes.
Then, in Section \ref{sec:disconch3} we investigate possible extensions of the setup of Sections \ref{sec:notatch3}-\ref{sec:pozch3} and we obtain, essentially as corollaries of the results of Section \ref{sec:mainch3}, explicit solutions to the Skorokhod embedding problem for spectrally one-sided L\'evy processes and their reflected versions. The last two sections are more technical in nature.
In Section \ref{sec:atomsch3} we discuss embeddings for measures with atoms, which are not covered by Theorem \ref{thm:mainch3}. Finally, in Section \ref{sec:proofsch3} we prove Theorems \ref{thm:mainch3} and \ref{thm:mainposch3}.

\section{Notation and basic setup}
\label{sec:notatch3}

We start by introducing the basic objects and notation that will be ubiquitous in this paper. We place ourselves in a general Markovian context and we follow closely Bertoin \cite{MR1406564} to which we refer for all the details. Specific notation connected with examples or particular cases will be introduced later, when necessary.

Let $(\Omega, \mathcal{F}, \p)$ be a probability space. We consider $X=(X_t:t\geq 0)$ a stochastic process taking values in some Polish space $(E,\rho)$ and having right-continuous sample paths. We want $X$ to be a 'nice' Markov process in the sense of Bertoin \cite{MR1406564}. To this end we denote $(\mathcal{F}_t:t\ge 0)$ its natural filtration taken right-continuous and completed. We suppose there is a family of probability measures $(\p_x:x\in E)$ on $\mathcal{F}$ such that for every stopping time $T<\infty$, under the conditional law $\p(\cdot|X_T=x)$, the shifted process $(X_{T+t}, t\geq 0)$ is independent of $\mathcal{F}_T$ and has the law of $X$ under $\p_x$. Furthermore we suppose $0$ is regular and instantaneous for $X$, that is $\p_0(\inf\{t>0:X_t=0\}=0)=1$ and $\p_0(\inf\{t\geq 0:X_t\neq 0\}=0)=1$, and that it is recurrent. We write $\p$ for $\p_0$ and $\e$ for $\e_0$, the expectation under $\p_0$. When we write $(X_t)$ we always mean the process in time $(X_t:t\ge 0)$. The process $B=(B_t:t\ge 0)$ denotes always a standard Brownian motion.

For a probability measure $\mu$ on $\re$, we denote its left-continuous tail by $\overline{\mu}(t)=\mu([t,\infty))$ and its support by $supp(\mu)$. The lower and upper bounds of the support are denoted respectively $a_\mu$ and $b_\mu$. Dirac's delta measure at point $y$ is denoted $\delta_y$.

All functions considered in the sequel are assumed to be Borel measurable.

\subsection{Markovian local time and excursions}
\label{sec:excch3}
We introduce now the local time and the excursion process of $X$ which will be our main tools in this paper. We follow Bertoin \cite{MR1406564} (see Blumenthal \cite{MR1138461} for an alternative approach based on potential theory).
The set of zeros of $X$ is denoted $\mathcal{Z}=\{t:X_t=0\}$. The last zero before time $t$ and the first zero after time $t$ are denoted respectively  $g_t=\sup\{u\leq t: X_u=0\}$, and $d_t=\inf\{u>t:X_u=0\}$. An interval of the form $(g_t,d_t)$ is called an excursion interval. These intervals appear in the canonical representation of the open set $[0,\infty)\setminus \overline{\mathcal{Z}}$ as the countable union of maximal disjoint open intervals.
The local time at $0$ of $X$ is denoted $L=(L_t:t\geq 0)$. Recall that it is characterized, up to a multiplicative constant, by the fact that the support of the Stieltjes measure $dL$ is included in $\overline{\mathcal{Z}}$ a.s.\ and that for any stopping time $T$, such that $X_T=0$ a.s.\ on $\{T<\infty\}$, the shifted process $(X_{T+t},L_{T+t}-L_T:t\geq 0)$, under $\p(\cdot|T<\infty)$, is independent of $\mathcal{F}_T$ and has the same law as $(X,L)$ under $\p$.
The local time $L_t$ is adapted to the filtration generated by zeros of $X$ (cf.\ Bertoin \cite[Thm.\ IV.4]{MR1406564}).
The right-continuous inverse of the local time $\tau=(\tau_l:l\geq 0)$, $\tau_l=\inf\{s\geq 0: L_s>l\}$, is a subordinator with infinite L\'evy measure $\Lambda^X$ (in particular $L_\infty=\infty$ a.s.). Note that the difference $(\tau_{l}-\tau_{l-})$ is just the length of the constancy period of $L$ at the level $l$, which in turn is the lifetime of the corresponding excursion. Thus $\overline{\mathcal{Z}}$ is the closure of the range of $(\tau_l:l\ge 0)$.

The space of excursions is defined as $U=\{\epsilon:\re_+\to E:\exists\; V(\epsilon),\; \epsilon(s)=0 \Leftrightarrow s\in \{0\}\cup [V(\epsilon),+\infty)\}$. The excursion process of $X$, $e=(e_l:l\geq 0)$ takes values in the space $U$ of excursions with an additional isolated point $\Upsilon$, that is $U\cup\{\Upsilon\}$, and is given by
\begin{equation}
  \label{eq:def_exc_prch3}
  e_l=(X_{\tau_{l-}+s},0\leq s\leq \tau_{l}-\tau_{l-})\quad \textrm{if}\quad \tau_{l-}<\tau_l,
\end{equation}
and $e_l=\Upsilon$ otherwise. One of the most important results for us, going back to the fundamental paper of It\^o \cite{MR0402949}, is that the above process is a Poisson point process with a certain characteristic measure $n$. This measure, called the It\^o measure, is uniquely determined up to a multiplicative constant factor. We will show however that our results are invariant under multiplication of the excursion measure by a constant.

The L\'evy measure of the subordinator $\tau$ can be easily deduced from It\^o's measure $n$. More precisely, as the lifetime of an excursion corresponds to the height of jump of $\tau$, we have $n(V(\epsilon)>a)=\Lambda^X((a,\infty))$, $a>0$ (we refer to Bertoin \cite[p.\ 117]{MR1406564} for details). Similar measures, with $V$ replaced by a general functional $F$, will be of prime importance in the sequel.

\subsection{Signed functionals of excursions}
\label{sec:funcch3}

We introduce now the main objects of this work, namely the class of signed functionals of Markovian excursions for which we want to solve the Skorokhod embedding problem.
We consider real, signed, continuous and monotone functionals of an excursion. An excursion is a function which starts in zero, travels in the space $E$ and then comes back to zero after some time (its lifetime). A generic excursion is taken constant and equal to zero after its lifetime. We stress that an excursion needs not to be continuous. Note that even when $E=\re$, an excursion can change sign during its lifetime.
The functionals we are interested in are just transformations of excursions. They are however real-valued, continuous and monotone. In particular, for a fixed excursion, they either stay positive or negative. At excursion's lifetime they take some value (the terminal value) and we define them to be constant and equal to the terminal value after excursion's lifetime.

More precisely let $F:U\times \re_+\to \re$. We will both write $F(\epsilon,t)$ and $F(\epsilon)(t)$, the latter being used to stress the time-dependence, with a particular excursion $\epsilon$ being fixed. For a fixed excursion $\epsilon\in U$, $F(\epsilon)$ is a \emph{monotone, continuous} function. It starts at zero, $F(\epsilon)(0)=0$, and is constant after the lifetime of $\epsilon$, that is $F(\epsilon)(t)=F(\epsilon)(V(\epsilon))$ for any $t\geq V(\epsilon)$.

Since we want the process induced by the functional to be adapted, we impose the condition that the value $F(\epsilon,t)$ is determined from the values of the excursion up to time $t$: $F(\epsilon,t)=F((\epsilon_s:s\leq t),t)$, that is for any $t\ge 0$ and $\epsilon,\epsilon'\in U$ such that $(\epsilon_s:s\leq t)=(\epsilon'_s:s\leq t)$, $F(\epsilon,t)=F(\epsilon',t)$. We set $F$ of the trivial excursion $\Upsilon$ to be zero: $F(\Upsilon)\equiv 0$, and assume that $|F(\epsilon,V(\epsilon))|>0$ for all $\epsilon\in U$. The last assumption can be relaxed as pointed out in Section \ref{sec:endch3}.

A functional $F$ induces a process in time $(F_t:t\ge 0)$, the value $F_t$ given as the functional $F$ of the excursion straddling time $t$ evaluated at the age of this excursion, that is
\begin{equation}
  \label{eq:def_F_tch3}
  \boxed{F_t=F\big(e_{L_t}\big)\big(t-g_t\big).}
\end{equation}
The assumption $|F(\epsilon, V(\epsilon))|>0$ implies that the set of zeros of the process $(F_t:t\ge 0)$ is equal to $\mathcal{Z}$, the set of zeros of $X$, and thus the local time $(L_t)$ is adapted to the natural filtration of $(F_t)$.
The process $(F_t:t\ge 0)$ is right-continuous with left limits and, when it jumps, it jumps to zero. This implies that for an excursion straddling time $s$, either $F_t$ is constant, and equal to its terminal value, for $t\in [t_0,d_s)$ for some $t_0$, $g_s\le t_0<d_s$, or $F_t$ converges to its terminal value $F_t\xrightarrow[t\to d_s]{a.s.} F(e_{L_s},V(e_{L_s}))$ but doesn't achieve it (since $F_{d_s}=0$).
The process of terminal values $\big(F(e_l,V(e_l)):l\geq 0\big)$ is a Poisson point process and we denote its characteristic measure by $n_F$. This measure is just the image of the excursion measure $n$ of $X$, by $\epsilon\to F(\epsilon,V(\epsilon))$.

Functional $F$ is called positive if $F:U\times \re_+\to\re_+$, that is if $F_t\ge 0$ for $t\ge 0$, or yet if $n_F((-\infty,0))=0$.

\emph{In this paper, when we speak of the process $F$ we always mean the process $(F_t:t\ge 0)$, denoted also $(F_t)$, which is adapted to the natural filtration of $X$. The natural filtration of $F$ designates the natural filtration of the process $(F_t:t\ge 0)$ that is $(\sigma\{F_s:s\le t\}:t\ge 0)$.
When speaking of characteristic measure of the functional associated with $(F_t)$ we mean the measure $n_F$.}

We close this subsection with some examples of functionals described above. Suppose that $X$ is a one-dimensional diffusion on $(l,r)$, $l<0<r$. In particular $X$ has continuous sample paths.
This implies that an excursion is either positive or negative and thus we can speak of $sgn(\epsilon)$, the sign of an excursion $\epsilon$. The first examples we present are connected with the age of excursion:
\begin{equation}
  \label{eq:sgnagech3}
  A(\epsilon)(t)=sgn(\epsilon)\big(t\land V(\epsilon)\big)\quad\textrm{and}\quad\alpha(\epsilon)(t)=\frac{\mathbf{1}_{sgn(\epsilon)=1}}{n_A((A(\epsilon,t),\infty))}-\frac{\mathbf{1}_{sgn(\epsilon)=-1}}{n_A((-\infty,A(\epsilon,t)))}.
\end{equation}
We use the notation $A(\cdot)$ and $\alpha(\cdot)$ because the two functionals are connected, one being a function of another. They yield the signed age process $A_t=sgn(X_t)(t-g_t)$ and the process
\begin{equation}
  \label{eq:genazemach3}
  \alpha_t=\frac{\mathbf{1}_{X_t>0}}{n_A((t-g_t,\infty))}-\frac{\mathbf{1}_{X_t<0}}{n_A((-\infty,g_t-t))}
\end{equation}
which is a martingale in the filtration generated by zeros of $X$, which is just the natural filtration of $(\alpha_t)$. Thus, the functional $\alpha$ yields a natural family of martingales $(\alpha_t)$ which are associated with the age process $(A_t)$.
That $(\alpha_t)$ is a martingale is well known for diffusions on natural scale (see Rainer \cite{MR1459486,MR1478743}, Pitman and Yor \cite[Rk.\ 3]{MR1478737}) but is generalized to our setup upon taking the scale function $s$ such that $s(0)=0$. Then $X$ and $s(X)$, which is a diffusion on the natural scale, have the same zeros and thus the same local times at zero and characteristic measures $n_A$ (up to a multiplicative constant).
The characteristic measure of $\alpha$ is easily seen to satisfy $n_\alpha(dv)=\frac{dv}{v^2}$. Note also that the measure $n_A$ satisfies $n_A((-\infty,-x)\cup (x,\infty))=n_V((x,\infty))=\Lambda^X((x,\infty))$, $x>0$. Furthermore, the measure $n_A$ is absolutely continuous with respect to the Lebesgue measure and its support is given as $\re$, $\re_+$ or $\re_-$ (cf.\ It\^o and McKean \cite[Sec.\ 6.2]{MR0345224}). In particular we can replace the open intervals in (\ref{eq:genazemach3}) with closed ones.
Profound studies of the It\^o measure $n$, in particular of the L\'evy measure $\Lambda^X$, were made via Krein's string theory. For more details on the measure $\Lambda^X$ we refer the reader to Knight \cite[pp.\ 71,77]{MR647781}, Kotani and Watanabe \cite{MR661628} and Bertoin \cite{MR1023954} (see Donati-Martin and Yor \cite{donati_yor_krein} for a recent account and further references).\\
When $X=B$ is a Brownian motion, then $\alpha_t=sgn(B_t)\sqrt{2\pi(t-g_t)}$, in which we recognize (up to a constant multiplicative factor) the celebrated Az\'ema martingale, which is the projection of $B$ on the filtration generated by its signs. Similar projection properties hold in the general setting (see Section \ref{sec:martch3}). Embeddings for these processes are discussed in detail in Section \ref{sec:azemach3}.\\
We can generalize upon (\ref{eq:sgnagech3}) in the abstract setting. Notice that even though it may not make any sense to speak of the sign of an excursion $\epsilon\in U$, the sign of the function $F(\epsilon)$ is well defined.
Thus, for a given functional $F$ we can define its scaled version $G^F$ through
\begin{equation}
  \label{eq:scaledfuncch3}
  G^F(\epsilon,t)=\frac{\mathbf{1}_{sgn(F(\epsilon))=1}}{n_F((F(\epsilon,t),\infty))}-\frac{\mathbf{1}_{sgn(F(\epsilon))=-1}}{n_F((-\infty,F(\epsilon,t)))}.
\end{equation}
When the measure $n_F$ is absolutely continuous with respect to the Lebesgue measure, it is immediate that $n_G(dx)=dx/x^2$, $x\neq 0$. We will come back to this matter in the remarks after Theorem \ref{thm:mainch3}.

A large family of functionals is given by
\begin{equation}
  \label{eq:intfuncch3}
 F^\gamma_\beta(\epsilon)(t)= sgn(\epsilon)\Big(\int_0^{t\land V(\epsilon)}|\epsilon(s)|^\beta ds\Big)^\gamma,
\end{equation}
where $\gamma,\beta$ are taken such that $F_t$ can be well defined. This family
is related with the Cauchy principal value associated with Brownian local times. We will refer to the functional $F^1_{-1}$, for $X=B$ a Brownian motion, as to the Brownian principal value. This is in fact an abuse of terminology as the associated process, at time $t$, is given via $\int_{g_t}^t \frac{ds}{B_s}$ which is an absolutely convergent integral. It is the process $P_t=\int_0^t\frac{ds}{B_s}=\lim_{x\to 0}\int_0^t \frac{ds}{B_s}\mathbf{1}_{|B_s|\ge x}$ which needs to be understood as Cauchy's principal value. However as the two processes are closely connected, the latter being the sum of the first one over excursions, we keep this naming convention.\\
The family given in (\ref{eq:intfuncch3}) contains also functional $F^1_1$ which is
connected with the area processes $\int_0^t |B_s|ds$ and $\int_0^t B^+_sds$, objects of great interest ever since the works of Cameron and Martin \cite{MR0013240} and Kac \cite{MR0016570}. We refer the reader to Perman and Wellner \cite{MR1422979} for a study using excursion theory and some applications in statistics.

The signed extrema process can be obtained taking
\begin{equation}
  \label{eq:extremach3}
  M(\epsilon)(t)=sgn(\epsilon)\sup_{s\le t\land V(\epsilon)}|\epsilon(s)|,
\end{equation}
which yields $M_t=sgn(X_t)\sup_{g_t\le s\le t}|X_s|$. This is a particularly important functional for us as our stopping times $T$, defined below in (\ref{eq:def_stopch3}), satisfy $X_T=M_T=\mathbf{1}_{X_T\ge 0}\sup_{s\le T}X_s + \mathbf{1}_{X_T<0}\inf_{s\le T}X_s$. Thus if we describe the distribution of $M_T$ we automatically describe the distribution of $X_T$.

When the underlying Markov process $X$ has the self-similarity property, $(X_{ct}:t\ge 0)\stackrel{\mathcal{L}}{=}(c^\kappa X_t:t\ge 0)$ for some $\kappa\in \re$, most of the functionals exemplified above fall into an important class of \emph{homogeneous functionals}. More precisely, following Carmona \textit{et al}.\ \cite{MR2001f:60080} and Pitman and Yor \cite{MR2002b:60139}, we say that $F$ is a $\theta$-\emph{homogeneous functional} of $X$ if
\begin{equation}
  \label{eq:homogfuncch3}
  F\Big(\epsilon, V(\epsilon)\Big)=V(\epsilon)^\theta F\Big(V(\epsilon)^{-\theta}\tilde{\epsilon}, 1\Big),
\end{equation}
where $\epsilon=(\epsilon_t:t\ge 0)\in U$ and $\tilde{\epsilon}=(\epsilon_{tV(\epsilon)}:t\ge 0)$, so that $V(\tilde{\epsilon})=1$.
The age and signed extrema are $\kappa$-homogeneous functionals and $F^\gamma_\beta$ is $\gamma(\beta\kappa+1)$-homogeneous. The characteristic measures of homogeneous functionals are easier to calculate thanks to the scaling property of $F$ inherited from $X$. This was exemplified in Ob\l\'oj and Yor \cite{obloj_yor}.

We close this section with some more remarks on the functionals and processes which can be treated in our setup. First, note that our study includes as well positive functionals. Numerous examples are obtained upon considering the absolute values of the functionals specified so far.\\
We pointed out above that with $X=B$, a Brownian motion, we can obtain in our setup the process $p_t=\int_{g_t}^t ds/B_s$ but not the process $P_t=\int_0^t ds/B_s$ (understood properly). However, we have to bear in mind that changing $X$ might allow us to treat such processes.
In particular we could consider $X_t=(R^{(q)}_t, K_t)$ where $R^{(q)}$ is a Bessel process with index $q\in (1,2)$ and $K$ appears in the classical Dirichlet process decomposition $R_t^{(q)}=B_t+\frac{(1-q)}{2} K_t$, and is locally of zero energy.
An independent definition of $K_t$ goes through the family of local times of $R^{(q)}$
and we can also write $K_t=\int_0^t ds/R^{(q)}_s$ where the integral can be taken as Cauchy principal value or \emph{partie finie} in Hadamard's sense associated with local times of $X$ (see Yor \cite[Sec.~10.4]{MR98e:60140} and Bertoin \cite{MR1030728} for details). Bertoin \cite{MR1030728} showed that $(0,0)$ is regular for $X$ and described It\^o's measure of excursions of $X$. As explained above, we could use the signed extrema functional (of the first coordinate of excursions of $X$ away from $(0,0)$) to control the process $(K_t)$. We note also that the process $K_t$ is actually a time-changed version of $\int_0^t ds/|B_s|^{1+1/q}$. We will not go further into this domain as it is not our aim here, but we hope these examples served to illustrate the generality of our setup.

\subsection{The Skorokhod embedding problem}
\label{sec:skoroch3}

The main aim of this paper is to solve the Skorokhod embedding problem for processes $(F_t)$ described above in Section \ref{sec:funcch3}. Let us recall the classical Skorokhod embedding problem as introduced by Skorokhod \cite{MR32:3082a} and developed by numerous authors afterwards.
 The problem is as follows: given a Brownian motion $(B_t:t\ge 0)$ and a probability measure $\mu$, find a stopping time $T$ in the natural filtration of $B$, such that $B_T\sim \mu$ and that $(B_{t\land T}:t\ge 0)$ is a uniformly integrable martingale. This is seen to be possible if and only if $\mu$ is centered. For all further details, as well as for an account of existing solutions, we refer the reader to our survey paper \cite{genealogia}.

Actually Skorokhod \cite{MR32:3082b} assumed that $\mu$ had finite second moment and required that $\e T<\infty$. However, authors soon realized that this was somehow artificial. It seems that the right way of saying that $T$ should be small (otherwise there is a trivial solution to the problem) is to require $T$ to be minimal. We recall that $T$ is minimal for process $(X_t)$ if for a stopping time $S$, $S\leq T$ and $X_S\sim X_T$ imply $S=T$ a.s. In the standard Brownian setup, minimality of $T$ is equivalent to the condition that $(B_{t\land T}:t\ge 0)$ is a uniformly integrable martingale (cf.\ Cox \cite{cox_min}, Cox and Hobson \cite{cox_hobson2}, Ob\l\'oj \cite[Sec.~8]{genealogia}).

Skorokhod \cite{MR32:3082b} originally considered randomized stopping times, which was sufficient for his needs. However, again, it was soon realized that it is much more natural to work with stopping times in the natural filtration of Brownian motion. More generally, when developing a solution to the Skorokhod embedding for some process $Y$ one should try to work with the stopping times in the natural filtration of $Y$.
Ob\l\'oj and Yor \cite{obloj_yor} developed a solution to the Skorokhod embedding for the age process of Brownian excursions. The study allowed them only to obtain a randomized solution for the embedding problem for the Az\'ema martingale. In this paper we improve upon their result by giving an explicit, non-randomized solution.

The problem we want to solve here is the following: for a given functional $F$, as in Section \ref{sec:funcch3}, and a probability measure $\mu$ on $\re$, describe explicitly a minimal stopping time $T$, in the natural filtration of $F$, such that $F_T\sim \mu$. The construction can require certain properties of the measure $\mu$ and we will say that $\mu$ is admissible if it has these properties.

Working within such a general setup may seem odd at first glance. Naturally, our study was motivated by examples, such as the Az\'ema martingale (see Section \ref{sec:azemach3} below). However, we have soon realized that our method applies very well in the general abstract setting described above and this actually allows to understand better the nature of our solution. In particular, as a special case, we will recover the solution obtained in Ob\l\'oj and Yor \cite{obloj_yor} and we will be able to explain the appearance, in their main theorem, of the \emph{dual Hardy-Littlewood} function (see (\ref{eq:dualHLch3}) below).

We will look for the solution to the Skorokhod embedding problem among the stopping times of the form already suggested in Ob\l\'oj and Yor \cite{obloj_yor}. Let $\varphi_-,\varphi_+:\re_+\to\re_+$ be two non-decreasing, right-continuous functions. Define
\begin{eqnarray}
  \label{eq:def_stopch3}
  T^F_{\varphi_-,\varphi_+}&=&T^F_{\varphi_-}\land T^F_{\varphi_+},\quad \textrm{where}\nonumber\\
  T^F_{\varphi_-}&=&\inf\{t>0:F_t\leq -\varphi_-(L_t)\},\\
  T^F_{\varphi_+}&=&\inf\{t>0:F_t\geq \varphi_+(L_t)\}.\nonumber
\end{eqnarray}
For a given probability measure $\mu$ on $\re$ we will search to specify the functions $\varphi_-$ and $\varphi_+$ such that $F_{T^F_{\varphi_-,\varphi_+}}\sim \mu$. We will write, when we want to stress a particular dependence, $T^F_{\varphi_-,\varphi_+}=T^F=T^F_\mu$. We will also drop the superscript $F$, when no confusion about the functional under consideration is possible.

We stress that $T^F$ is a stopping time in the natural filtration of $(F_t)$. This follows from the fact that zeros of $(F_t)$ are equal to zeros of $X$ and thus the local time $(L_t)$ is adapted to the natural filtration of $(F_t)$. This in turn was a consequence of the assumption $|F(\epsilon,V(\epsilon))|>0$ imposed in Section \ref{sec:funcch3}. Everything that follows can be easily extended to the case when this assumption is removed, but at the cost of considering stopping times in the filtration of $(F_t,L_t)$ or of $(X_t)$.

\section{Main results}
\label{sec:mainch3}

In this Section we present our main theorems. Theorem \ref{thm:mainch3} gives an explicit solution to the Skorokhod embedding problem for signed functionals and non-atomic probability measures. Propositions \ref{cor:mainagech3} and \ref{prop:embdiffch3} are applications of Theorem \ref{thm:mainch3} respectively for the signed age process $(A_t)$ of a diffusion, and for a diffusion on natural scale. They contain necessary and sufficient conditions for the existence of an embedding which are similar to the classical phrasing of the Skorokhod embedding problem. In Section \ref{sec:azemach3} we present more applications of Theorem \ref{thm:mainch3}, we obtain solutions to the Skorokhod embedding problem for Az\'ema's martingale, for Brownian principal value, for Brownian motion itself and for its two-sided extrema process.

Theorem \ref{thm:mainposch3} deals with positive functionals and presents an explicit solution for an arbitrary probability measure on $\re_+$. Naturally the choice of working with positive rather then negative functionals is arbitrary and the results can be easily translated for negative functionals. In Section \ref{sec:pozch3} we will apply Theorem \ref{thm:mainposch3} to describe explicit solutions to the Skorokhod embedding problems for Bessel processes of dimension $\delta\in (0,2)$, for their maximum processes and for age processes of excursions. We will also cover the age process of excursions of Cox-Ingersoll-Ross processes.

\subsection{Signed functionals}
\label{sec:mainsignedch3}

With a non-atomic probability measure $\mu$ on $\re$ and a functional $F$, as described in Section \ref{sec:funcch3}, we associate the following functions:
\begin{equation}
  \label{eq:int_much3}
  D_{\mu}(y)=\int_{[0,y]} \frac{d\mu(s)}{n_F([s,+\infty))}\quad\textrm{and}\quad G_{\mu}(x)=\int_{[x,0]}\frac{d\mu(s)}{n_F((-\infty,s])},
\end{equation}
for $y\ge 0$ and $x\le 0$, and where $n_F$ is the characteristic measure of Poisson point process of the terminal values of $F$, as defined in Section \ref{sec:funcch3}.
The inverses $D^{-1}_\mu$, $G^{-1}_\mu$ are taken right-continuous, $D^{-1}_\mu|_{[D_\mu(b_\mu),\infty)}=\infty$, $G^{-1}_\mu|_{[G_\mu(a_\mu),\infty)}=-\infty$. We make the following fundamental assumptions
\begin{eqnarray}
\label{eq:suppFch3}
  x,y\in supp(\mu)\Rightarrow n_F((-\infty,x])\cdot n_F([y,\infty))>0,\\
 D_\mu(\infty)=\int_0^\infty \frac{d\mu(s)}{n_F([s,+\infty))}=\int_{-\infty}^0 \frac{d\mu(s)}{n_F((-\infty,s])}=G_\mu(-\infty),  \label{eq:cond_regch3}
\end{eqnarray}
which ensures that the functions $D^{-1}_\mu(G_\mu(\cdot))$ and $G^{-1}_\mu(D_\mu(\cdot))$ are well defined. Thus, for $y,z\ge 0$, we can define
\begin{eqnarray}
  \label{eq:psi_expressch3}
  \psi_+(y)&=&
\int_0^y\frac{d\mu(s)}{n_F\Big([s,+\infty)\Big)\Big(1+\overline{\mu}(s)-\overline{\mu}\big(G^{-1}_{\mu}(D_\mu(s))\big)\Big)}\\
\psi_-(z)&=&\int^0_{{\scriptscriptstyle -}z}\frac{d\mu(s)}{n_F\Big((-\infty,s]\Big)\Big(1+\overline{\mu}\big(D_\mu^{-1}(G_\mu(s))\big)-\overline{\mu}(s)\Big)}\label{eq:psi2_expressch3}
\end{eqnarray}
and ${\psi_-}_{|_{\re_-}}\equiv {\psi_+}_{|_{\re_-}}\equiv 0$. As a consequence of (\ref{eq:finiteproofch3}) below, we will see that
$\psi_+(y)=\psi_-(z)=\infty$ for $y\ge b_\mu$ and $z\ge -a_\mu$. Define the left-continuous inverses $\varphi_-,\varphi_+:\re_+\to\re_+$ by
\begin{eqnarray}
  \label{eq:varphisch3}
\varphi_+(y)&:=&\psi_+^{-1}(y)=\inf\{x\ge 0: \psi_+(x)\ge y\}\nonumber\\
\varphi_-(y)&:=&\psi_-^{-1}(y)=\inf\{x\ge 0: \psi_-(x)\ge y\},
\end{eqnarray}
$\varphi_-(0)=\varphi_+(0)=0$. When there will be two or more functionals considered, we will often add a superscript $F$ to the functions defined in (\ref{eq:int_much3}-\ref{eq:varphisch3}) to avoid any confusion.\\
We are ready to present our main result.
\begin{theorem}
\label{thm:mainch3}
  Let $F$ be a functional as defined in Section \ref{sec:funcch3} and $\mu$ a non-atomic probability measure on $\re$ such that $\mu(\re_-)>0$, $\mu(\re_+)>0$,
and (\ref{eq:suppFch3}) and (\ref{eq:cond_regch3}) hold. Then
\begin{equation}
  \label{eq:def_stop1ch3}
  \begin{split}
   T^F_{\varphi_-,\varphi_+} & =\inf\Big\{t>0:\psi_-(-F_t)\ge L_t\textrm{ or }\psi_+(F_t)\ge L_t\Big\}\\
&=\inf\Big\{t>0:F_t\notin (-\varphi_-(L_t),\varphi_+(L_t))\Big\},
  \end{split}
\end{equation}
where $\psi_{+/-}$, $\varphi_{+/-}$ are given by (\ref{eq:psi_expressch3})-(\ref{eq:varphisch3}), is an a.s.\ finite stopping time and it solves the Skorokhod embedding problem for $(F_t:t\ge 0)$, i.e.\ $F_{T^F_{\varphi_-,\varphi_+}}\sim \mu$.
Furthermore, $T^F=T^F_{\varphi_-,\varphi_+}$ is minimal and $F_{T^F}=\sup_{t\le T^F}F_t\cdot \mathbf{1}_{F_{T^F}\ge 0}+\inf_{t\le T^F}F_t\cdot \mathbf{1}_{F_{T^F}<0}$.
\end{theorem}
\noindent \textbf{Some remarks about Theorem \ref{thm:mainch3}}

\emph{The characteristic measure $n$ and the local time $L$ are defined up to a multiplicative constant. However our solution is invariant under a renormalization of $n_F$ and $L$ as proved in Section \ref{sec:solch3}.}

\emph{The equality between two expressions for $T^F_{\varphi_-,\varphi_+}$ in (\ref{eq:def_stop1ch3}) follows from the fact that $\mu$ is non-atomic. The first form in somewhat more explicit, however the second one is, in a sense, more universal. When we treat measures with atoms, as in Theorem \ref{thm:mainposch3} below, only the analogue of the second form is valid, hence the notation $T^F_{\varphi_-,\varphi_+}$ above or $T^F_\varphi$ in (\ref{eq:stopposch3}). Taking $\varphi_{+/-}$ left continuous is consistent with weak inequalities ($\le$) in (\ref{eq:def_stop1ch3}) and left-continuous tails of the measure $n_F$ in (\ref{eq:int_much3})-(\ref{eq:psi2_expressch3}). Our results can be naturally re-written in the ``right-continuous'' convention. We note that taking $\varphi_{+/-}$ right-continuous doesn't affect (a.s.) $T^F_{\varphi_-,\varphi_+}$ since the law of $L_{T^F_{\varphi_-,\varphi_+}}$ is absolutely continuous with respect to the Lebesgue measure, see (\ref{eq:law_L_Tch3}).}

\emph{The solution described in Theorem \ref{thm:mainch3} depends on $F$ only through its characteristic measure $n_F$. We note that in the particular case when $n_F(dx)=dx/x^2$ the formulae (\ref{eq:int_much3}-\ref{eq:psi2_expressch3}) simplify considerably.
Recall that to a functional $F$ we associated its scaled version $G^F$ through (\ref{eq:scaledfuncch3}) and that $n_G(dx)=dx/x^2$. In a sense then we can always work with the particular case of the characteristic measure $dx/x^2$ if we agree to solve the Skorokhod embedding for $G^F$ instead of $F$.}

\emph{As one would expect, in the symmetric case the expressions in Theorem \ref{thm:mainch3} simplify significantly. More precisely, suppose that $n_F$ and $\mu$ are symmetric, i.e.\ invariant under $x\to-x$. Then the stopping time in (\ref{eq:def_stop1ch3}) can be written as $T^F=\inf\{t>0: \psi(|F_t|)\ge L_t\}$, where $\psi(y)=\int_0^y\frac{d\mu(s)}{2\overline{\mu}(s)n_F([s,\infty))}$. This yields a simple justification of the solution for positive functionals presented in Theorem \ref{thm:mainposch3} in the non-atomic case.}

\emph{The assumption, in Theorem \ref{thm:mainch3}, that the measure $\mu$ has no atoms is important. We will treat the case of measures with atoms separately in Section \ref{sec:atomsch3} and we will see that the formulae become more involved.}

\emph{We stress the property that the stopping time $T^F_{\varphi_-,\varphi_+}$ is minimal. We recall that this is the general requirement imposed on a solution to the Skorokhod embedding problem (see Section \ref{sec:skoroch3} above or Ob\l\'oj \cite[Sec.~8]{genealogia} for all the details). It generalizes the traditionally imposed condition of uniform integrability of a certain stopped martingale.}

\emph{The solution in the special case of Brownian motion and the signed extrema functional (\ref{eq:extremach3}) yields a solution to the Skorokhod embedding problem for Brownian motion which coincides with the solution of Vallois \cite{MR86m:60200}. See Section \ref{sec:martch3} below for details.
This solution to the Skorokhod embedding problem for Brownian motion $B$ is such that the terminal value $B_T$ is equal either to the maximum $\sup_{s\leq T}B_s$ or to the minimum $\inf_{s\le T}B_s$. This is reminiscent of the solution developed by Perkins \cite{MR88k:60085}. We note however that our embedding relies actually on a third process, namely the local time L, while Perkins' solution is expressed solely in terms of Brownian motion, its maximum and its minimum.}

\emph{The functions $\psi_{+/-}$ are taken increasing. Using the same methodology, we could also develop an analogue embedding but with $\psi_{+/-}$ decreasing. In Brownian motion setup this was done by Vallois \cite{MR1162722} and we will come back to this issue after Proposition \ref{prop:bound_ltagech3}.}
\smallskip

We stress that the solution presented in Theorem \ref{thm:mainch3} is very general and requires only the knowledge of the characteristic measure $n_F$. In Section \ref{sec:azemach3} below we will study the special case of $n_F(dx)=\frac{dx}{2x^2}\mathbf{1}_{x\ne 0}$. This is the characteristic measure of a number of functionals including the functionals associated with Az\'ema's martingale $\tilde{\alpha}_t=\sqrt{\frac{\pi}{2}}sgn(B_t)\sqrt{t-g_t}$ and the signed maximum $M_t=sgn(B_t)\sup_{g_t\le s\le t}|B_s|$ process.

We specify now to the functionals related to the age process of excursions. In this case we dispose of a family of martingales (\ref{eq:genazemach3}) which allows us to understand better the condition (\ref{eq:cond_regch3}).
\begin{proposition}
\label{cor:mainagech3}
Let $(X_t:t\ge 0)$ be a one-dimensional diffusion on $(l,r)$, $l<0<r$, $X_0=0$ a.s., and $A_t=sgn(X_t)(t-g_t)$ be the signed age process of excursions of $X$. Recall that its scaled version $(\alpha_t:t\ge 0)$, given in (\ref{eq:genazemach3}), is a martingale.
For a non-atomic probability measure $\mu$ on $\re$ there exists a stopping time $S$ in the natural filtration of $(A_t)$ such that $A_S\sim \mu$ and $(\alpha_{t\land S}:t\ge 0)$ is a uniformly integrable martingale \emph{if and only if} $D^A_\mu(\infty)=G^A_\mu(-\infty)<\infty$.\\
If $D^A_\mu(\infty)=G^A_\mu(-\infty)<\infty$ then $T^A_{\varphi_-,\varphi_+}$, given in (\ref{eq:def_stop1ch3}), is an a.s.\ finite stopping time in the natural filtration of $(A_t)$, $A_{T^A_{\varphi_-,\varphi_+}}\sim \mu$, and $(\alpha_{t\land T^A_{\varphi_-,\varphi_+}}:t\ge 0)$ is a uniformly integrable martingale.
\end{proposition}
\noindent The proof of Proposition \ref{cor:mainagech3} is presented in Section \ref{sec:azemach3}. We just note here that the condition $\e \alpha_{T^A_{\varphi_-,\varphi_+}}=0$ is equivalent to $D^A_\mu(\infty)=G_\mu^A(-\infty)<\infty$, which is a version of (\ref{eq:cond_regch3}).
We recalled in Section \ref{sec:skoroch3} that in the standard Skorokhod embedding problem the condition that the stopping time $S$ should be \emph{small} is imposed by requiring that a certain martingale, stopped at $S$, should be \emph{uniformly integrable}. To $(A_t:t\ge 0)$, the process of age of excursions, a natural family of martingales $(\alpha_t:t\ge 0)$, displayed in (\ref{eq:genazemach3}), is associated. Proposition \ref{cor:mainagech3} shows that our criterion (\ref{eq:cond_regch3}) for the age process $A$, agrees with the condition of uniform integrability of stopped martingale $\alpha$. 

Typically, one obtains a solution to the Skorokhod embedding problem for diffusions, by adapting a solution developed for Brownian motion (cf.~Ob\l\'oj \cite[Sec.~8]{genealogia}). Here we obtain an embedding for diffusions directly from Theorem \ref{thm:mainch3}. This demonstrates an advantage of the abstract phrasing of Theorem \ref{thm:mainch3}.
\begin{proposition}
\label{prop:embdiffch3}
  Let $(X_t:t\ge 0)$ be a one-dimensional diffusion on $(l,r)$, $l<0<r$, $X_0=0$ a.s. Assume that $X$ is on natural scale and choose the classical normalization of the local time $L$ given by $L_t=|X_t|-\int_0^t sgn(X_s)dX_s$ a.s. Let $M$ be the signed extrema functional given by (\ref{eq:extremach3}), $\mu$ a non-atomic centered probability measure on $(l,r)$, and $\psi_{+/-}$, $\varphi_{+/-}$ defined via (\ref{eq:psi_expressch3})-(\ref{eq:varphisch3}).
Then the characteristic measure of $M$ is given by $n_M(dx)=\frac{dx}{2x^2}$, $x\neq 0$.
The stopping time $T^M_{\varphi_-,\varphi_+}$, defined in (\ref{eq:def_stop1ch3}), satisfies
  \begin{eqnarray}
    \label{eq:stop_diffch3}
    T^M_{\varphi_-,\varphi_+}&=& \inf\Big\{t>0:\psi_-(-X_s)\geq L_t\textrm{ or }\psi_+(X_s)\ge L_t\Big\}\nonumber\\
&=&\inf\Big\{t>0:X_t\notin (-\varphi_-(L_t),\varphi_+(L_t))\Big\},
  \end{eqnarray}
$X_{T^M_{\varphi_-,\varphi_+}}\sim\mu$, and $(X_{t\land T^M_{\varphi_-,\varphi_+}}:t\ge 0)$ is a uniformly integrable martingale.
\end{proposition}
\noindent The condition that $\mu$ is centered, that is $\int_{\re} |x|d\mu(x)<\infty$ and $\int_{\re} xd\mu(x)=0$, is a necessary condition for existence of a stopping time $S$ such that $X_S\sim\mu$ and $(X_{t\land S}:t\ge 0)$ is a uniformly integrable martingale (cf.~\cite[Sec.~8]{genealogia}). It is equivalent to $D_\mu^M(\infty)=G^M_\mu(-\infty)<\infty$.

The proof of Proposition \ref{prop:embdiffch3} is given is Section \ref{sec:azemach3}. Essentially we have to prove that $n_M(dx)=dx/2x^2$ and then apply Theorem \ref{thm:mainch3}. Note that the special form of $T^M_{\varphi_-,\varphi_+}$ implies that it induces the same embedding both for $(M_t)$ and for $(X_t)$.

\subsection{Positive functionals}
\label{sec:mainposch3}

In Ob\l\'oj and Yor \cite{obloj_yor}, which inspired the present study, the authors considered positive functionals of Brownian motion with some particular interest placed upon the functionals $F$ with $n_F(dx)=\frac{dx}{x^2}\mathbf{1}_{x> 0}$.
We will now extend this. We will see that in this case we can have an explicit formulae for arbitrary measures. We will use a similar methodology as Ob\l\'oj and Yor \cite{obloj_yor} and it comes as no surprise that we recover their results.
However, we think that the general setup studied in this paper allows to understand better the particular formulae obtained in \cite{obloj_yor}.
In Section \ref{sec:pozch3} we will develop applications for functionals $F$ with $n_F(dx)=\frac{dx}{x^2}\mathbf{1}_{x> 0}$ and provide in (\ref{eq:dualHLch3}) a corrected form of the \emph{dual Hardy-Littlewood} function introduced in \cite{obloj_yor}.
\begin{theorem}
\label{thm:mainposch3}
  Let $F$ be a positive functional, as defined in Section \ref{sec:funcch3}, and $\mu$ a probability measure on $\re_+$ with $\mu(\{0\})=0$ and $n_F([y,\infty))>0$ for $y\in supp(\mu)$.
Define
  \begin{equation}
    \label{eq:psiposch3}
\psi_\mu(y)=\int_0^y\frac{\mathbf{1}_{\mu(\{s\})=0}d\mu(s)}{\overline{\mu}(s)n_F([s,\infty))}+\sum_{s< y}\frac{\ln\Big(\frac{\overline{\mu}(s)}{\overline{\mu}(s+)}\Big)}{n_F([s,\infty))}\mathbf{1}_{\mu(\{s\})>0},
  \end{equation}
and $\varphi_\mu$ its right-continuous inverse. Then the stopping time
\begin{equation}
  \label{eq:stoppos1ch3}
  T^F_{\varphi_\mu}=\inf\big\{t>0:F_t\ge \varphi_\mu(L_t)\big\}
\end{equation}
is a.s.\ finite and solves the Skorokhod embedding problem for $F$, i.e.\ $F_{T^F_{\varphi_\mu}}\sim\mu$. Furthermore, $T^F_{\varphi_\mu}$ is minimal and $\sup_{t\le T^F_{\varphi_\mu}}F_t=F_{T^F_{\varphi_\mu}}$.\\
For $\mu$ a probability measure on $\re_+$ with $\mu(\{0\})=\varsigma>0$ define $\tilde{\mu}=\mu-\varsigma (\delta_{0}-\delta_{\infty})$. Then $\psi_{\tilde{\mu}}(\infty)<\infty$ and the stopping time
\begin{equation}
  \label{eq:stopposch3}
  \tilde{T}^F_{\mu}=T^F_{\varphi_{\tilde{\mu}}}\land \inf\big\{t>0:L_t=\psi_{\tilde{\mu}}(\infty)\big\}
\end{equation}
embeds $\mu$ in $F$, i.e.\ $F_{\tilde{T}^F_\mu}\sim \mu$.
\end{theorem}
Taking above $\varphi_\mu$ right-continuous is more convenient for the proof. However, as for the Theorem \ref{thm:mainch3}, we note that since the law of $L_{T^F_{\varphi_\mu}}$ is proved to be absolutely continuous with respect to the Lebesgue measure, taking $\varphi_\mu$ left-continuous does not affect our solution.

The second part of the theorem provides a way of dealing with an atom at zero of $\mu$ which follows the idea of Vallois \cite{MR86m:60200}. Note that if we applied the first part of the theorem for $\mu$ with $\mu(\{0\})>0$, we would have $T^F_{\varphi_\mu}=0$ a.s.
Another way of dealing with an atom at zero is to use a standard external randomization (cf.\ Ob\l\'oj \cite[Sec.~6.1]{genealogia}).

Note that like in Theorem \ref{thm:mainposch3} the stopping times we define are minimal which is the property required from a solution to the Skorokhod embedding problem (cf.~\cite[Sec.~8]{genealogia}).

We will see that Theorem \ref{thm:mainposch3} is significantly simpler to prove than Theorem \ref{thm:mainch3}. Actually, upon taking $G_\mu\equiv 0$, for probability measure $\mu$ with $\mu(\re_-)=0$, Theorem \ref{thm:mainposch3} for non-atomic measures can be seen as a direct corollary of Theorem \ref{thm:mainch3}. In particular, the expression (\ref{eq:psiposch3}) for $\psi_\mu$ then coincides with the expression (\ref{eq:psi_expressch3}) for $\psi_+$.

Finally we stress that although we chose to work with positive functionals, one could just as well work with negative functionals. All our results for positive functionals have an obvious rewriting in terms of negative functionals.

\subsection{Links with martingale theory}
\label{sec:martch3}

We come back to the last remarks below Theorem \ref{thm:mainch3} and establish a link with the solution to the Skorokhod embedding problem developed by Vallois \cite{MR86m:60200}. More generally,
as we will rely on excursion theory throughout the paper, we mention possible martingale theory arguments and establish a link with the study of Jeulin and Yor \cite{MR83c:60110} which was the cornerstone of the work of Vallois \cite{MR86m:60200}. We recall that Ob\l\'oj and Yor \cite{obloj_yor} presented martingale theory arguments and already established there the link with Vallois' solution for symmetric measures. Here we complete this discussion.

In Proposition \ref{prop:embdiffch3} we showed how to apply Theorem \ref{thm:mainch3} to obtain an embedding for a one-dimensional diffusion. Consider then the particular case of Brownian motion. The stopping time in (\ref{eq:stop_diffch3}) can be written as $T_{\varphi_-,\varphi_+}=\inf\{t>0: B_t=-\varphi_-(L_t)\textrm{ or }B_t=\varphi_+(L_t)\}$. We recognize instantly the form of the stopping times considered by Jeulin and Yor \cite{MR83c:60110} and exploited in Vallois \cite{MR86m:60200}. Recall that as we work with the signed extrema functional which obeys $n_M(dx)=dx/2x^2$ we have $D_\mu(y)=2\int_0^y sd\mu(s)$ and $G_\mu(x)=-2\int_x^0sd\mu(s)$. These coincide with the functions $2\rho^+$ and $2\rho^-$ defined in Vallois \cite{MR86m:60200}. Suppose that $\mu$ has a positive density so that the functions $D_\mu$, $G_\mu$, $\psi_{+/-}$ are continuous and strictly increasing. It is then a matter of simple calculation to check from the nested definition of functions $h^+$ and $h^-$ in \cite{MR86m:60200} that their inverses satisfy the same differential equations as $\psi_{+/-}$ and thus to see, a posteriori, that our solution coincides with the solution of Vallois.

Our presentation of the embedding for Brownian motion has the advantage of being more explicit than the one developed by Vallois \cite{MR86m:60200}. His solution, on the other hand, works for arbitrary centered probability measures on $\re$, whereas our formulation breaks down in the presence of atoms, as explained in Section \ref{sec:atomsch3} below.

As mentioned above, Vallois \cite{MR86m:60200} exploited the work of Jeulin and Yor \cite{MR83c:60110} who described the law of the couple $(B_{T_{\varphi_-,\varphi_+}},T_{\varphi_-,\varphi_+})$ for a large class of functions $\varphi_{+/-}$. Jeulin and Yor \cite{MR83c:60110} used stochastic calculus, much in the spirit of Az\'ema and Yor \cite{MR82c:60073a}, and special families of Brownian martingales. Recently, Nikeghbali \cite{ashkan_class} generalized some results of Az\'ema \cite{MR88f:60081} and Jeulin and Yor \cite{MR83c:60110} and was able to obtain a solution to the Skorokhod embedding problem for Bessel processes of dimension $\delta\in (0,2)$. We recover his results in Proposition \ref{cor:bessch3}, which is a corollary of Theorem \ref{thm:mainposch3}.\\
Our present work provides a parallel to Jeulin and Yor \cite{MR83c:60110} and Nikeghbali \cite{ashkan_class} replacing martingale theory arguments with excursion theory arguments. Such a possibility was already mentioned in Jeulin and Yor \cite{MR83c:60110}. However, it seems to us that using excursion theory allows us to widen significantly the scope of the investigation and develop a more general solution to the Skorokhod embedding problem as compared with the results one can obtain following the original approach of Jeulin and Yor \cite{MR83c:60110}.

To end this section, we want to present some remarkable martingales and apply them to establish optimal properties of our stopping times.
Consider $(X_t)$ a one-dimensional diffusion, martingale on some open interval containing zero and let $(\mathcal{H}_t)$ be the right-continuous version of $(\mathcal{F}_{g_t})$.
Rainer \cite{MR1459486} showed that the optional projection of $X^+_t$ on $\mathcal{H}_t$ is given as $\frac{1}{2}\alpha^+_t$, where $(\alpha_t)$ is the martingale in (\ref{eq:genazemach3}). Assume furthermore that $n_A(dv)$ is \emph{symmetric}. Recall that $K(L_t)-|B_t|k(L_t)$ is a martingale for any bounded Borel function $k$, where $K(y)=\int_0^y k(s)ds$ and $(B_t)$ is a Brownian motion and $(L_t)$ is its local time at zero (see Ob\l\'oj \cite{ja_martcarac}). This generalizes instantly, through Dambis-Dubins-Schwarz theorem, to $X$ in place of $B$. Thus for any bounded Borel function $k$, the process $M^k_t=K(L_t)-|X_t|k(L_t)$ is a martingale, when now $L_t$ is the local time of $X$ normalized so $(|X_t|-L_t)$ is a martingale. Since the local time $(L_t)$ is adapted to $(\mathcal{H}_t)$, the projection of $(M^k_t)$ on $(\mathcal{H}_t)$ is given as
\begin{equation}
  \label{eq:mart_proj_agech3}
m^k_t=K(L_t)-\frac{k(L_t)}{2n_A([t-g_t,\infty))}=K(L_t)-\frac{k(L_t)|\alpha_t|}{2},\quad t\ge 0,
\end{equation}
 and is a $(\mathcal{H}_t)$-martingale. Ob\l\'oj and Yor \cite{ja_yor_maxmart} used martingales $M^k$ with $k(x)=\mathbf{1}_{x\ge \lambda}$, $\lambda\ge 0$, to establish a bound on the law of $L_R$ when the law of $|X_R|$ is fixed. We can use the same argumentation, but with martingales $m^k$ instead of $M^k$, to obtain a parallel result, this time with the distribution of $|\alpha_R|$ fixed.
 \begin{proposition}
\label{prop:bound_ltagech3}
   Let $X$ be a one-dimensional diffusion, martingale on $(l,r)$, $l<0<r$ such that $n_A(dv)$ is symmetric. Let $\mu$ be a probability measure on $(0,\infty)$ with $\int_0^\infty xd\mu(x)<\infty$ and $R$ be a stopping time such that $|\alpha_R|\sim \mu$ and $(\alpha_{t\land R}:t\ge 0)$ is a uniformly integrable martingale. Denote $\rho_R$ the law of $L_R$. The following bound is true
\begin{equation}
  \label{eq:bound_ltch3}
\e\Big[\Big(L_R-\overline{\rho}_R^{-1}(p)\Big)^+\Big]\leq \e\Big[\Big(L_{T^{|\alpha|}}-\overline{\rho}_{T^{|\alpha|}}^{-1}(p^*)\Big)^+\Big],\quad p\in [0,1],
\end{equation}
where $T^{|\alpha|}=T^{|\alpha|}_{\varphi_\mu}$ is given in Theorem \ref{thm:mainposch3}, the inverses $\overline{\rho}^{-1}_\cdot$ are taken left-continuous and  $p^*=\overline{\mu}\Big(\overline{\mu}^{-1}(p)\Big)\ge p$.
 \end{proposition}
\noindent The proof of Proposition \ref{prop:bound_ltagech3} is just an application of the optional stopping theorem to the martingale in (\ref{eq:mart_proj_agech3}) with $k(l)=\mathbf{1}_{l\ge \overline{\rho}_R^{-1}(p)}$. We have
\begin{equation*}\label{eq:proof_prop_age}
    \begin{split}
\e\Big[\Big(L_R-\overline{\rho}_R^{-1}(p)\Big)^+\Big]
&=\frac{1}{2}\e|\alpha_R|\mathbf{1}_{L_R\ge \overline{\rho}_R^{-1}(p)}\le \frac{1}{2}\e|\alpha_R|\mathbf{1}_{|\alpha_R|\ge \overline{\mu}^{-1}(p)}\\
&=\frac{1}{2}\e|\alpha_{T^{|\alpha|}}|\mathbf{1}_{|\alpha_{T^{|\alpha|}}|\ge \overline{\mu}^{-1}(p)}
\le\frac{1}{2}\e|\alpha_{T^{|\alpha|}}|\mathbf{1}_{L_{T^{|\alpha|}}\ge \overline{\rho}_{T^{|\alpha|}}^{-1}(p^*)}\\
&=\e\Big[\Big(L_{T^{|\alpha|}}-\overline{\rho}_{T^{|\alpha|}}^{-1}(p^*)\Big)^+\Big],
\end{split}
\end{equation*}
which proves the Proposition. Note that in the statement we could just as well fix the law of $|A_R|$ as it is equivalent to fixing the law of $|\alpha_R|=1/n_A([|A_R|,\infty))$, only the integrability condition on $\mu$ would change.

Proposition \ref{prop:bound_ltagech3} and the bound of Ob\l\'oj and Yor \cite{ja_yor_maxmart} can be summarized by saying that the law of $L_\infty$, the local time of a continuous UI martingale $(N_t)$, the distribution of $|N_\infty|$ or $|A_\infty|$ being fixed, is bounded in the \emph{excess wealth order} and hence in the \emph{convex order} (see Kochar \textit{et al.} \cite{MR2003i:60024}) and the upper bound is attained with stopping times of the type (\ref{eq:stopposch3}). This complements the study of Vallois \cite{MR1162722} who obtained lower and upper bounds on the law of $L_\infty$ in the \emph{convex order}, under fixed distribution of $N_\infty$ and showed that both bounds can be attained with solutions to the Skorokhod embedding problem for Brownian motion. The upper bound is attained with the solution developed by Vallois \cite{MR86m:60200} which we recovered in Theorem \ref{thm:mainch3}. The lower bound was attained with an analogous solution presented by Vallois \cite{MR1162722} which, in comparison with (\ref{eq:def_stop1ch3}), takes the functions $\psi_{+/-}$ decreasing and not increasing. As noted after Theorem \ref{thm:mainch3}, we could also develop our general solution with decreasing functions $\psi_{+/-}$. This would lead to lower bounds on the law of $L_\infty$ under fixed law of $|A_\infty|$.

\section{Applications to  the Az\'ema martingale and other signed functionals}
\label{sec:azemach3}

In this section we develop applications of Theorem \ref{thm:mainch3} to various functionals with specific characteristic measures. We deal in particular with the Az\'ema martingale, which actually motivated our study. We also prove Propositions \ref{cor:mainagech3} and \ref{prop:embdiffch3}. We use the notation from Section \ref{sec:mainch3}.

Let $X$ be a one-dimensional diffusion on $(l,r)$, $l<0<r$. Recall the martingale $(\alpha_t:t\ge 0)$ displayed in (\ref{eq:genazemach3}) and the fact that the characteristic measure of the functional $\alpha$ is given by $n_\alpha(dv)=\frac{dv}{v^2}$, $v\neq 0$. The following proposition is, in a sense, a variant of Proposition \ref{cor:mainagech3}.
\begin{proposition}
\label{cor:agemartch3}
  For a non-atomic measure $\mu$ on $\re$ there exists a stopping time $S$ such that $\alpha_S\sim \mu$ and $(\alpha_{t\land S}:t\ge 0)$ is a uniformly integrable martingale \emph{if and only if} the measure $\mu$ is centered, in which case we can take $S=T^\alpha_{\varphi_-,\varphi_+}$ defined via (\ref{eq:def_stop1ch3}) and (\ref{eq:psi_expressch3}-\ref{eq:varphisch3}).
\end{proposition}
\textbf{Proof. }
We prove Proposition \ref{cor:mainagech3}. Proposition \ref{cor:agemartch3} then follows.
Suppose $S$ is a stopping time such that $A_S\sim\mu$ and $(\alpha_{t\land S}:t\ge 0)$ is a uniformly integrable martingale. Then  $\e |\alpha_S|<\infty$ and $\e \alpha_S=0$. We have
  \begin{eqnarray}
    \label{eq:helppropproofch3}
  D_\mu^A(y)=\int_0^y \frac{d\mu(v)}{n_A([v,\infty))}=
\e \frac{A_S\mathbf{1}_{0\le A_S\le y}}{n_A((A_S,\infty))}=\e \alpha_S\mathbf{1}_{0\le A_S\le y},
  \end{eqnarray}
where the second equality follows from the fact that $A_S\sim \mu$ and $n_A$ is absolutely continuous with respect to the Lebesgue measure (cf.\ Section \ref{sec:funcch3}). In parallel with (\ref{eq:helppropproofch3}), we obtain $G_\mu^A(x)=\e\alpha_S\mathbf{1}_{x\le A_S<0}$. We see that $\e \alpha_S=0$ is equivalent to $D_\mu^A(\infty)=G^A_\mu(-\infty)<\infty$.
 However then the condition (\ref{eq:cond_regch3}) is satisfied and Theorem \ref{thm:mainch3} tells us that $T=T^A_{\varphi_-^A,\varphi_+^A}<\infty$ a.s.\ and $A_T\sim \mu$.
Let $T_n$ denote $T^A_{\varphi_-^A\land n,\varphi_+^A\land n}$. The process $(\alpha_{t\land T_n}:t\ge 0)$ is bounded and hence a uniformly integrable martingale. Furthermore, $T_n\to T$ so a sufficient condition for the uniform integrability of $(\alpha_{t\land T}:t\ge 0)$ is the uniform integrability of $(|\alpha_{T_n}|:n\ge 1)$. This in turn will follow from $L^1$ convergence of the sequence. We have $|\alpha_{T_n}|\to|\alpha_T|$ a.s. By Sheffe's lemma, it suffices to show $\e |\alpha_{T_n}|\to\e |\alpha_T|$, which follows from
\begin{eqnarray*}
\e|\alpha_{T_n}|&=&\e|\alpha_T|\mathbf{1}_{|\alpha_T|<n}+n\p(|\alpha_T|\ge n)\\&=&\int_{-n}^n |x|d\mu(x)+n\int_{|x|\ge n}d\mu(x)\le \int_{\re}|x|d\mu(x)=\e|\alpha_T|
\end{eqnarray*}
Thus, $(\alpha_{t\land T}:t\ge 0)$ is a uniformly integrable martingale if and only if $\alpha_T$ is integrable (and thus centered), which in turn is equivalent to $D_\mu^A(\infty)=G^A_\mu(-\infty)<\infty$. This ends the proof of Proposition \ref{cor:mainagech3}.
The proof of Proposition \ref{cor:agemartch3} is similar, it suffices to note that as $n_\alpha(dx)=\frac{dx}{x^2}\mathbf{1}_{x\neq 0}$, the condition $D_\mu^\alpha(\infty)=G_\mu^\alpha(-\infty)<\infty$ is equivalent to $\int_{\re}|x|d\mu(x)<\infty$ and $\int_{\re}xd\mu(x)=0$. The uniform integrability of $(\alpha_{t\land T}:t\ge 0)$ is argued as above.

Likewise, the proof of Proposition \ref{prop:embdiffch3} is immediate if we show that $n_M(dy)=\frac{dy}{2y^2}$. There are number of ways to see it. Here is one. From Corollary 2.1 and discussion in Section 2.4 in Pitman and Yor \cite{MR1439532} it follows that $n_M((y,\infty))$ is proportional to the scale function $s(y)$ of the diffusion $(Y_t:t\ge 0)$ obtained upon conditioning $X$, starting at some point $x>0$, to approach $\infty$ before $0$. More precisely, we define $Y$ as Doob's $h$-transform of $X$ via $\e_x[H(Y_s:s\leq t)]=\e_x[\frac{X_{t\land T_0}}{x}H(X_s:s\leq t)]$, where $H$ is a positive functional and $T_r=\inf\{t:X_t=r\}$. Taking $H(X_s:s\leq t)=\mathbf{1}_{t\land T_b<T_a}$, $0<a<x<b$, using the fact that $X$ is on natural scale and letting $t\to\infty$, it is easy to verify that $s(y)= 1/y$ is a scale function for $Y$. Similar argument applies with $(-X)$ in place of $X$, which shows that $n_M$ is symmetric. Thus $n_M(dx)$ is proportional to $dx/x^2$, $n_M(dx)=cdx/x^2$ for some positive constant $c$. \\
Note that in Proposition \ref{prop:embdiffch3} we chose a specific normalization for the local time, described by the fact that $(L_t-|X_t|)$ is a local martingale.
This allows us to recover $c$ using the compensation formula (cf.\ Revuz and Yor \cite[Prop.~XII.2.6]{MR2000h:60050}). Let $\lambda>0$ and $T_{-\lambda,\lambda}=\inf\{t:|X_t|\ge \lambda\}$, then we obtain
$$1=\e\Big[\sum_{g_s\le T_{-\lambda,\lambda}}\mathbf{1}_{\sup_{g_s\le u\le d_s} X_u\ge \lambda}+\mathbf{1}_{\inf_{g_s\le u\le d_s} X_u\le -\lambda}\Big]=2c\frac{\e L_{T_{-\lambda,\lambda}}}{\lambda}=2c,$$
which yields $c=\frac{1}{2}$ and ends the proof.
$\Box$\smallskip

We now specialize to the Brownian setup.
Let $B=(B_t:t\ge 0)$ be a real-valued Brownian motion and $L=(L_t:t\ge 0)$ its local time at zero with $\e L_t=\e|B_t|$. Define functional $p$ through $p(\epsilon,t)=\int_0^{t\land V(\epsilon)}\frac{ds}{\epsilon(s)}$. The process $p_t$ is just equal to $\int_{g_t}^t\frac{ds}{B_s}$, and is closely linked with the Cauchy principal value associated with Brownian local times $P_t=\int_0^t\frac{ds}{B_s}$, where the integral is understood as the principal value (see Section \ref{sec:funcch3} above, cf.\ Biane and Yor \cite{MR88g:60188}). Set $\tilde{p}=\frac{1}{2}p_t$. Introduce also $\tilde{\alpha}_t=\sqrt{\frac{\pi}{2}}sgn(B_t)\sqrt{t-g_t}$ the Az\'ema martingale and recall that it is the projection of $B$ on the filtration generated by its zeros.
\begin{proposition}
\label{cor:mainch3}
  Let $\mu$ be a non-atomic probability measure on $\re$, such that $\int_{\re_+} |s|d\mu(s)=\int_{\re_-} |s|d\mu(s)$. Define $\psi_{+/-}$, $\varphi_{+/-}$ via (\ref{eq:psi_expressch3})-(\ref{eq:varphisch3}).
For a process $(F_t)$, define the stopping time
  \begin{equation}
    \label{eq:defstopspecch3}
    \begin{split}
    T(F)&=\inf\Big\{t>0:\psi_-(-F_t)\ge L_t\textrm{ or }\psi_+(F_t)\ge L_t\Big\}\\&=
\inf\Big\{t>0: F_t\notin (-\varphi_-(L_t),\varphi_+(L_t))\Big\}.
    \end{split}
  \end{equation}
Then we have $B_{T(B)}\sim \mu$, $\tilde{p}_{T(\tilde{p})}\sim\mu$ and $\tilde{\alpha}_{T(\tilde{\alpha})}\sim \mu$. Furthermore, $T(B)$, $T(\tilde{p})$ and $T(\tilde{\alpha})$ are stopping times in the natural filtration of $B$, $\tilde{p}$ and $\tilde{\alpha}$ respectively. The martingales $(B_{t\land T(B)}:t\ge 0)$ and $(\tilde{\alpha}_{t\land T(\tilde{\alpha})}:t\ge 0)$ are uniformly integrable if and only if $\int_{\re}|x|d\mu(x)<\infty$, in which case $\mu$ is centered.
\end{proposition}
\textbf{Proof. }
We have $X=B$ a Brownian motion. Recall the notation of Section \ref{sec:mainch3}. As $\varphi_+$ and $\varphi_-$ are increasing, we have $M_{T^F_{\varphi_-,\varphi_+}}=B_{T^F_{\varphi_-,\varphi_+}}$, where $M_t=sgn(B_t)\sup_{g_t\le s\le t}|B_s|$ is associated with the functional $M$ given  by (\ref{eq:extremach3}).
We will now show that
\begin{equation}
 \label{eq:charmeas3ch3}
 n_M(dx)=n_{\tilde{\alpha}}(dx)=n_{\tilde{p}}(dx)=\frac{dx}{x^2}\mathbf{1}_{x\neq 0}.
\end{equation}
This will end the proof, as then we can proceed exactly as in the proof of Propositions \ref{cor:mainagech3} and \ref{cor:agemartch3} above. Here we can also exploit the well known fact that for $\e|B_T|<\infty$ the conditions: $T$ is minimal and $(B_{t\land T}:t\ge 0)$ is a uniformly integrable martingale, are equivalent (cf.~Ob\l\'oj \cite[Sec.~8]{genealogia}).

The assertion for $n_{\tilde{\alpha}}$ is a direct consequence of the independence between the length and the sign of an excursion and the expression of the  characteristic measure of the age functional, $V(\epsilon)(t)=t\land V(\epsilon)$, given by $n_V=\frac{dv}{\sqrt{2\pi v^3}}$ (cf.\ Revuz and Yor \cite[Prop.\ XII.2.8]{MR2000h:60050}). It follows also from (\ref{eq:bessmaxch3}) and (\ref{eq:bessviech3}) with $q=-\frac{1}{2}$.
Note that $\alpha_t=sgn(B_t)\sqrt{2\pi(t-g_t)}=2\tilde{\alpha}_t$, where $\alpha_t$ is defined via (\ref{eq:genazemach3}). \\
The assertion (\ref{eq:charmeas3ch3}) for $n_M$ is well known (cf.\ Revuz and Yor \cite[Prop.\ XII.3.6]{MR2000h:60050}) and we argued it in a greater generality above in the proof of Proposition \ref{cor:agemartch3}.\\
Finally, the assertion on the characteristic measure of $\tilde{p}$ follows readily from Theorem 4.1 in Biane and Yor \cite{MR88g:60188}, but we present another simple justification. We look at the process of terminal values of $p$, $p_{\tau_l-}=\int_0^{\tau_l-\tau_{l-}}\frac{ds}{B_{\tau_{l-}+s}}$. It follows easily from the Poisson point process properties of the excursion process, that $H_l=P_{\tau_l}=\sum_{u\leq l}p_{\tau_u-}$ is a L\'evy process. Examining the scaling property for $H$ one finds that $\frac{1}{\pi}H_l$ is actually a standard Cauchy process. Now using the exponential formula we can calculate $\e \mathrm{e}^{\mathrm{i}\lambda p_t}$ and comparing it with the known quantity for Cauchy process, we find $n_p(dx)=\frac{dx}{x^2}\mathbf{1}_{x\ne 0}$ and thus $n_{\tilde{p}}(dx)=\frac{dx}{2x^2}$, $x\neq 0$.
$\Box$\smallskip

Note that the above Proposition is quite remarkable as we have actually the same formula for the stopping time which works both for Brownian motion and for its projection on the filtration generated by the signs that is the Az\'ema martingale.

The last process we want to mention here is the skew Brownian motion. Intuitively speaking, it is a Brownian motion which chooses positive excursions with probability $p$, $0<p<1$, and negative excursions with probability $(1-p)$. Naturally, we know the characteristic measure of the signed extrema functional (\ref{eq:extremach3}), $n_M(dx)=\frac{pdx}{x^2}\mathbf{1}_{x>0}+\frac{(1-p)dx}{x^2}\mathbf{1}_{x<0}$, and we can thus develop an explicit solution to the Skorokhod embedding problem for skewed Brownian motion.\\
Skew Brownian motion, or more general skewed processes, were used recently by Cox and Hobson \cite{cox_hobson_skew} to develop a class of embeddings in Brownian motion. Their solution contains both the solutions of Vallois \cite{MR86m:60200} and of Az\'ema and Yor \cite{MR82c:60073a}.
More precisely, Cox and Hobson \cite{cox_hobson_skew} exploit general skew processes, solutions of $X_t=B_t-G(L^X_t)$, where $G$ is a function with Lipschitz constant at most one. In light of their work, it might be interesting to consider embeddings in such processes, however when $G'$ is not a constant such processes do not have Markov property and they can not be treated with our methodology.

We close this section with two explicit calculations of functions $\varphi_{+/-}$ for different probability measures $\mu$. We restrain ourselves to the case $n_F(dx)=dx/x^2$, $x\neq 0$, which we encountered in Propositions \ref{prop:embdiffch3}, \ref{cor:agemartch3} and \ref{cor:mainch3}. We look only on asymmetric probability measures as the symmetric case follows immediately from the solution for positive functionals which are given in Section \ref{sec:pozch3} below.
\smallskip\\
\textbf{Double exponential.} Let $\mu(dx)=\frac{\lambda^2}{\lambda+\gamma}\mathrm{e}^{-\lambda x}\mathbf{1}_{x>0}+\frac{\gamma^2}{\lambda+\gamma}\mathrm{e}^{\gamma x}\mathbf{1}_{x<0}$, for some $\lambda,\gamma>0$. The coefficients are chosen so that $\mu$ is a centered probability measure. We have $D_\mu(y)=(1-\mathrm{e}^{-\lambda y}(1+\lambda y))/(\lambda+\gamma)$, $y\ge 0$, and $G_\mu(x)=(1-\mathrm{e}^{\gamma x}(1-\gamma x))/(\lambda+\gamma)$, $x\le 0$. Note that $D_\mu(\infty)=1/(\lambda+\gamma)=G_\mu(-\infty)$ so that (\ref{eq:cond_regch3}) is indeed verified. We see easily that $G_\mu(\lambda x/\gamma)=D_\mu(-x)$ and $D_\mu(\gamma y/\lambda)=G_\mu(-y)$. This yields
  \begin{equation}
    \label{eq:regexampexpch3}
    \varphi_-(y)=\frac{1}{\lambda}\sqrt{(\lambda+\gamma)y}\quad\textrm{and}\quad
\varphi_+(y)=\frac{1}{\gamma}\sqrt{(\lambda+\gamma)y}.
  \end{equation}
\noindent\textbf{F-uniform.} Let $\mu(dx)=K(\mathbf{1}_{x\ge g}/x^2+\mathbf{1}_{x\le -h}/x^2)dx$ where $g,h>0$ and $K(1/g+1/h)=1$. We called this measure \emph{F-uniform} as it is just a (weighted) restriction of the measure $n_F(dx)$ to $\re\setminus (-h,g)$. In particular it is easy to justify that $\varphi_{+/-}$ have to be affine functions.
We recover the formulae presented in Vallois \cite{MR86m:60200}. We have $D_\mu(y)=K\log(y/g)$, $y\ge g$, and $G_\mu(x)=K\log(-x/h)$, $x<-h$. Note that $G_\mu(-\infty)=D_\mu(\infty)=\infty$ so (\ref{eq:cond_regch3}) is verified. We obtain easily
\begin{equation}
  \label{eq:regexamppowerch3}
  \varphi_+(y)=g\Big(1+\frac{y}{2K}\Big)\quad\textrm{and}\quad \varphi_-(y)=h\Big(1+\frac{y}{2K}\Big).
\end{equation}

\section{Application to Bessel and Cox-Ingersoll-Ross processes}
\label{sec:pozch3}
In this section , we apply Theorem \ref{thm:mainposch3} to obtain embeddings for positive functionals. More precisely, we specialize to the case of $F$ with $n_F(dx)=\frac{dx}{x^2}\mathbf{1}_{x>0}$, which was studied in Ob\l\'oj and Yor \cite{obloj_yor}. We have $n_F([y,\infty))=1/y$ and so we obtain
\begin{equation}
\label{eq:dualHLch3}
\psi_\mu(t)=\int_0^t\frac{s\mathbf{1}_{\mu(\{s\})=0}d\mu(s)}{\overline{\mu}(s)}+\sum_{s\le y}s\ln\Big(\frac{\overline{\mu}(s)}{\overline{\mu}(s+)}\Big)\mathbf{1}_{\mu(\{s\})>0},
  \end{equation}
which is the correct definition of the \emph{dual Hardy-Littlewood function}, introduced in Ob\l\'oj and Yor \cite[Eq.~3.1]{obloj_yor}.
We now understand well the appearance of the identity function $s$ under the integral. In fact, it is a direct consequence of taking $n_F(dx)=\frac{dx}{x^2}\mathbf{1}_{x> 0}$, which is equivalent to saying $1/n_F([s,\infty))=s$, $s\ge 0$.

We will apply the general solution presented in Theorem \ref{thm:mainposch3} to obtain a solution to the Skorokhod embedding problem for Bessel processes and some of their functionals. The results we present were discovered independently by Nikeghbali \cite{ashkan_class}, who generalized work of Ob\l\'oj and Yor \cite{obloj_yor} using martingale theory and general theory of processes.

Let $(R^{(q)}_t:t\ge 0)$ be a Bessel process with index $q\in (-1,0)$, starting in zero (we write BES$^{(q)}$). For background on Bessel processes we refer to Revuz and Yor \cite[Ch.\ XI]{MR2000h:60050}.
We recall that zero is an instantly reflecting barrier for $R^{(q)}$ and that the so-called dimension of $R^{(q)}$ is given through $\delta=2q+2\in (0,2)$. The processes $R^{(q)}$ for $q\in (-1,0)$ are nice Markov processes as in Section \ref{sec:notatch3}, for which we have a convenient description of their It\^o's measures $n^q$, and thus we can use our embedding described in Theorem \ref{thm:mainposch3}. As $R^{(q)}$ is positive it is natural to redefine the space of excursion $U$ as the space of positive excursions. We recall that $M(\epsilon)$ signifies the maximum of an excursion and $V(\epsilon)$ its lifetime. The characteristic measures of these two functionals are denoted respectively $n^q_M$ and $n^q_V$.

The following two equivalent descriptions of $n^q$, due to Pitman and Yor \cite{MR656509} (see also Biane and Yor \cite[p.\ 43]{MR88g:60188}), provide an important generalization of Williams' decomposition of Brownian excursions (cf.\ Williams \cite[Sec.\ 67]{MR531031}, Rogers \cite{MR84h:60140}). In the sequel, we choose the normalization of It\^o's measure and the local time under which the process $(R^{(q)}_t)^{-2q}-L_t$ is a martingale (see Donati-Martin \emph{et al}.\ \cite{donati_roy_val_yor_norm} for a survey of common normalizations). We recall however that our results are independent of renormalization of the local time and It\^o's measure.
\\
\textbf{The first description of $n^q$, $q\in (-1,0)$.}
\begin{itemize}
\item The characteristic measure $n^q_M$ satisfies
  \begin{equation}
    \label{eq:bessmaxch3}
    n^q_M([x,\infty))=x^{2q},\quad x\ge 0,
  \end{equation}
\item for any $x>0$, conditionally on $M=x$, the maximum $M$ is attained in a unique instant $S$, $0<S<V$ a.s., and $(\epsilon(s):s\le S)$ and $(\epsilon(V-s):s\le V-S)$ are two independent BES$^{(-q)}$ processes stopped at their first hitting times of the level $x$.
\end{itemize}
\textbf{The second description of $n^q$, $q\in (-1,0)$.}
\begin{itemize}
\item The characteristic measure $n^q_V$ satisfies
  \begin{equation}
    \label{eq:bessviech3}
    n^q_V([v,\infty))=\frac{2^qv^q}{\Gamma(|q|+1)},\quad v\ge 0,
  \end{equation}
\item for any $v>0$, conditionally on $V=v$, the process $(\epsilon(t):t\le v)$ is a Bessel bridge with index $(-q)$, going from $0$ to $0$ on time interval of length $v$.
\end{itemize}
Consider two functionals of the excursions of $R^{(q)}$, given via
\begin{equation}
  \label{eq:funcbessch3}
\tilde{M}(\epsilon)(t)=\Big(\sup_{s\le t\land V(\epsilon)}\epsilon(s)\Big)^{2|q|}\quad \textrm{and}\quad \tilde{A}(\epsilon)(t)=\Big(t\land V(\epsilon)\Big)^{|q|},
\end{equation}
and note that their characteristic measures are given via $n_{\tilde{M}}(dx)=\frac{dx}{x^2}\mathbf{1}_{x>0}$ and $n_{\tilde{A}}(dv)=\frac{cdv}{x^2}\mathbf{1}_{v>0}$, where $c=\frac{2^q}{\Gamma(|q|+1)}$.

We can apply Theorem \ref{thm:mainposch3} generalizing
the results of Ob\l\'oj and Yor \cite{obloj_yor} from Brownian motion to any Bessel process with index $q\in (-1,0)$. We note that the same result was obtained independently by Nikeghbali \cite{ashkan_class} using entirely different methods. For simplicity, we treat below the case of measures without atom in zero but, similarly to Theorem \ref{thm:mainposch3}, this is not necessary.
\begin{proposition}
  \label{cor:bessch3}
Let $\mu$ be a probability measure on $\re_+$ with $\mu(\{0\})=0$. Define the \emph{dual Hardy-Littlewood} function $\psi_\mu$ through (\ref{eq:dualHLch3}) and let $\varphi_\mu$ denote its right-continuous inverse. Then the stopping times
\begin{eqnarray}
  \label{eq:bessstopch3}
  T^R&=&\inf\big\{t>0:\tilde{M}_t\ge \varphi_\mu(L_t)\big\}=\inf\big\{t>0:R^{(q)}_t\ge \varphi_\mu(L_t)^{1/2|q|}\big\}\nonumber \\
&=&\inf\big\{t>0:\sup_{s\le t} R^{(q)}_s \ge \varphi_\mu(L_t)^{1/2|q|}\big\},\\
 T^{\tilde{A}}&=&\inf\big\{t>0:\tilde{A}_t\ge \varphi_\mu(cL_t)\big\}=\inf\big\{t>0:(t-g_t)\ge \varphi_\mu(cL_t)^{1/|q|}\big\},\nonumber
\end{eqnarray}
solve the Skorokhod embedding problem, i.e.\ $\tilde{M}_{T^R}=\big(\sup_{s\le T^R}R^{(q)}_s\big)^{2|q|}=\big(R^{(q)}_{T^R}\big)^{2|q|}\sim\mu$ and $\tilde{A}_{T^{\tilde{A}}}\sim \mu$, where $q\in (-1,0)$ and $c=2^q/\Gamma(|q|+1)$.
\end{proposition}
\textbf{Proof. }
  The embedding for $\tilde{M}$ is a simple application of Theorem \ref{thm:mainposch3}. As $\varphi_\mu$ is increasing it is straightforward to see that $\tilde{M}_{T_R}=\big(\sup_{s\le T^R}R^{(q)}_s\big)^{2|q|},$ and as the local time $L_t$ is constant on excursions away from zero, $T^R$ is a point of increase for the maximum process of $R^{(q)}$ and thus $R^{(q)}_{T_R}=\sup_{s\le T^R}R^{(q)}_s$.\\
To prove the embedding for $\tilde{A}$ is suffices to notice that the function $\psi_\mu^{\tilde{A}}$ associated to $\tilde{A}$ by (\ref{eq:psiposch3}) is linked with the \emph{dual Hardy Littlewood} function (\ref{eq:dualHLch3}) through $\psi_\mu^{\tilde{A}}=\frac{1}{c}\psi_\mu$, and thus their inverses satisfy $\varphi^{\tilde{A}}_\mu(y)=\varphi_\mu(cy)$.
$\Box$\smallskip

The discussion above can be extended to Bessel processes with drift downwards. More precisely, Donati-Martin and Yor \cite{donati_yor_krein} showed that the measure
\begin{equation}
  \label{eq:bessdriftch3}
n^{q,\alpha}(dv)=C\frac{\exp(-\beta v)}{v^{1-q}}dv,\quad -1<q\le 0,\; \beta>0,
\end{equation}
on $\re_+$ can be seen as the L\'evy measure of the inverse of the local time at zero of BES$(q,\beta\downarrow)$ process, that is downwards BES$^{(-q)}$ process with ``drift'' $\beta$. We refer to Pitman and Yor \cite{MR620995} for definition of BES$(q,\beta\downarrow)$ processes. Thus (\ref{eq:bessdriftch3}) gives the characteristic measure of the age functional of excursions of BES$(q,\beta\downarrow)$ process and we can apply Theorem \ref{thm:mainposch3} to obtain an explicit solution to the Skorokhod embedding problem for the age process $A_t=t-g_t$ for $X$ a BES$(q,\beta\downarrow)$ process.

We turn now to an embedding for the age process of excursions for the Cox-Ingersoll-Ross processes. Fix $\gamma>0$, $\delta\in(0,2)$  and let $(X_t:t\ge 0)$ be the non-negative solution of
\begin{equation}
  \label{eq:circh3}
  dX_t=2\sqrt{X_t}dB_t+(\delta-2\gamma X_t)dt,
\end{equation}
where we assume $X_0=0$, and where $(B_t)$ is a real-valued Brownian motion.
The processes $X=X^{\delta,\gamma}$ found applications in mathematical finance (cf.\ Yor \cite[Ch.\ 5]{MR1854494}) and are called the \emph{Cox-Ingersoll-Ross} processes. They are also called the \emph{squared Ornstein-Uhlenbeck process} with dimension $\delta$ (cf.\ Pitman and Yor \cite{MR1478737}).
 Denote $\Lambda=\Lambda^{\delta, \gamma}$  the L\'evy measure of $\tau$, the inverse of the local time $L$ at zero of $X$. We recall that $\Lambda^{\delta,\gamma}$ is also the characteristic measure of the age functional $A$ given in (\ref{eq:sgnagech3}), $\Lambda((v,\infty))=n_A((v,\infty))$. This measure is known and given by
 \begin{equation}
   \label{eq:measurecirch3}
   \Lambda^{\delta, \gamma}((v,\infty))=C\frac{\mathrm{e}^{-2\gamma(1-\frac{\delta}{2}) v}}{\big(1-\mathrm{e}^{-2\gamma v}\big)^{(1-\frac{\delta}{2})}},
 \end{equation}
where $C$ is a constant which depends on the normalization of the local time $L$ (cf.\ Pitman and Yor \cite[Eq.\ 59]{MR1478737}). More precisely multiplying $L$ by $c$, divides $n$, and therefore also $\Lambda$, by $c$. As proved in Section \ref{sec:solch3}, our solution to the Skorokhod embedding problem is independent of such renormalizations. A possible canonical choice of $C$ is given by $C=2\gamma \big(\Gamma(\delta/2)\Gamma(1-\delta/2)\big)^{-1}$, see Pitman and Yor \cite[Sec.\ 4]{MR1478737} for the details.

Applying Theorem \ref{thm:mainposch3}, we obtain therefore instantly an explicit solution to the Skorokhod embedding problem for age process $A_t=t-g_t$ of the Cox-Ingersoll-Ross process.

We close this section with some examples of measures $\mu$ and the corresponding \emph{dual Hardy-Littlewood} functions. These examples are found in Ob\l\'oj and Yor \cite{obloj_yor} and are reported here for the sake of completeness. However, the formula for the geometrical law given in \cite{obloj_yor} was incorrect and we display here its corrected form.\\
\textbf{Weibull's law.} Take $a>0$, $b\ge 0$ and let $\mu(dx)=abx^{b-1}\mathrm{e}^{-a x^b}\mathbf{1}_{[0,\infty)}(x)dx$, so that the tail is equal to $\overline{\mu}(x)=\mathrm{e}^{-ax^b}$. Simple calculations show that
$$\psi_\mu(x)=\frac{ab}{b+1}x^{b+1},\quad \varphi_\mu(x)=\Big(\frac{b+1}{ab}x\Big)^{\frac{1}{b+1}}.$$
In particular, for $b=0$, $\mu$ is the exponential law with parameter $a$ and we have $\phi_\mu(x)=\sqrt{2x/a}$.\\
\textbf{Uniform law.} Let $\mu$ be the uniform law on $[a,b]$. We have
$$\psi_\mu(x)=\int_a^x\frac{y}{b-y}dy=\Big(b\log\Big(\frac{b-a}{b-x}\Big)-(x-a)\Big),$$ for $a\leq x< b$ and we put $\psi(x)=0$ for $x\in [0,a)$ and $\psi(x)=\infty$ for $x\geq b$.
The inverse function is not given by any explicit formula, as it would be equivalent to solving $x=1-c\mathrm{e}^{-x}$, where $c$ is a certain constant. However, as the target measure has no atoms, we can still write the stopping time as:
$$T_\mu=\inf\Big\{t> 0:\psi_\mu\big(F_t\big)\geq L_t\Big\},$$
where $F_t$ is the functional with $n_F(dx)=\frac{dx}{x^2}\mathbf{1}_{x>0}$ for which we develop the embedding.\\
\textbf{Geometric law}. Let $\mu$ be a probability measure on $\nr$ with $\mu(\{k\})=(1-p)^{k-1}p$, for certain $0<p<1$, $k\in\nr$. Then $\overline{\mu}(k)=(1-p)^{k-1}$ and
\begin{eqnarray*}
&&\psi_\mu(k)=\frac{k(k+1)}{2}\log\Big(\frac{1}{1-p}\Big)\\
&&\varphi_\mu(x)=k\quad \textrm{for}\quad\frac{(k-1)k}{2}\log\Big(\frac{1}{1-p}\Big)\leq x<\frac{k(k+1)}{2}\log\Big(\frac{1}{1-p}\Big).
\end{eqnarray*}
More generally, for any discrete probability measure we have $\varphi_\mu(x)=k$ for $x\in[a_{k-1},a_k)$, where $a_k=\psi_\mu(k)$.

\section{Discontinuous functionals: spectrally negative L\'evy processes}
\label{sec:disconch3}
In this section we widen the scope of the results presented so far.
We extend previously obtained embeddings to some cases when the underlaying process has jumps and thus the extrema functional (\ref{eq:extremach3}) is discontinuous. Our main goal is to give an explicit solution to Skorokhod embedding for spectrally negative L\'evy processes and their reflected versions.
The basic idea, coming from Pistorius \cite{pistorius_skoro}, is that when using embeddings presented in Theorems \ref{thm:mainch3} and \ref{thm:mainposch3} it is just important to preserve the distribution of the stopped local time and the equality that the stopped process is a function of the local time (with appropriate signs in the signed case). In this way we can use this embedding for some process $(X_t)$ even if its extrema functional has jumps, if only we can ensure that the process $(X_t)$ will come back to a given level, before hitting zero. We will now phrase this concept precisely, first in the case of positive functionals and reflected L\'evy processes and then in the signed case for spectrally negative L\'evy processes.\\
Let $(Y_t)$ be a $[0,\infty)$-valued nice Markov process (as described in Section \ref{sec:notatch3}) and $M$ its maximum functional given in (\ref{eq:extremach3}), so that $M_t=\sup_{g_t\le s\le t}Y_s$. Assume furthermore that the following property holds:
\begin{equation}
\label{eq:discassumch3}
    M_{t-}<m<M_t\Rightarrow \exists\,s: t\le s\le d_t,\; Y_s=m,
\end{equation}
where $d_t=\inf\{u>t: Y_u=0\}$. We then have the following corollary of Theorem \ref{thm:mainposch3}.
\begin{corollary}
\label{cor:discgenposch3}
Let $(Y_t)$ be a $[0,\infty)$-valued Markov process as in Section \ref{sec:notatch3} which verifies (\ref{eq:discassumch3}) and $\mu$ a probability measure on $(0,\infty)$ with $n_M([y,\infty))>0$ for $y\in supp(\mu)$. Let $\varphi_\mu$ be the right-continuous inverse of $\psi_\mu$ defined in (\ref{eq:psiposch3}).
Then the stopping time
\begin{equation}\label{eq:stop_discch3}
    \overline{T}_{\varphi_\mu}=\inf\big\{t>0:Y_t = \varphi_\mu(L_t)\big\}
\end{equation}
is finite a.s.\ and solves the Skorokhod embedding problem for $(Y_t)$, i.e.\ $Y_{\overline{T}_{\varphi_\mu}}\sim \mu$.
\end{corollary}
The corollary follows from Theorem \ref{thm:mainposch3} and its proof presented in Section \ref{sec:proofposch3}. It suffices to note that $L_{\overline{T}_{\varphi_\mu}}= L_{T^M_{\varphi_\mu}}$, where $T^M_{\varphi_\mu}$ is defined in (\ref{eq:stoppos1ch3}), and that $Y_{\overline{T}_{\varphi_\mu}}=\varphi_\mu(L_{\overline{T}_{\varphi_\mu}})$. An extension of the above corollary to the case of measures $\mu$ with an atom in zero is also immediate along the lines of Theorem \ref{thm:mainposch3}.

We discuss now the case of $Y$ which is a reflected spectrally negative L\'evy process. Let $(X_t)$ be a spectrally negative L\'evy process which does not drift to $-\infty$ and $Y_t=\sup_{s\le t}X_s-X_t$. Define the scale function $W:[0,\infty)\to [0,\infty)$ via its Laplace transform
\begin{equation}\label{eq:scalefuncch3}
    \int_0^\infty \mathrm{e}^{-\theta x}W(x)dx=\frac{t}{\log \e[\exp(\theta X_t)]},
\end{equation}
which is well defined for $\theta$ large enough\footnote{We refer to Bertoin \cite{MR1406564} or Pistorius \cite{pistorius_skoro} for details.}.
Recall that the left- and right- hand derivatives $W_-'$ and $W_+'$ of the scale function are well defined and that the characteristic measure of the maximum functional is given as $n_M([y,\infty))=W_+'(y)/W(y)$ (cf.\ Pistorius \cite{MR2054585}). We can now rephrase Corollary \ref{cor:discgenposch3} for reflected spectrally negative L\'evy process.
\begin{proposition}\label{prop:refllevych3}
Let $(X_t)$ be a spectrally negative L\'evy process which does not drift to $-\infty$ and $Y_t=\sup_{s\le t}X_s-X_t$. For a probability measure $\mu$ on $(0,\infty)$ define
  \begin{equation}
    \label{eq:psiposdiscch3}
\psi_\mu(y)=\int_0^y\frac{\mathbf{1}_{\mu(\{s\})=0}W(s)}{\overline{\mu}(s)W_+'(s)}d\mu(s)
+
\sum_{s<y}\frac{\ln\Big(\frac{\overline{\mu}(s)}{\overline{\mu}(s+)}\Big)W(s)}{W_+'(s)}\mathbf{1}_{\mu(\{s\})>0},
  \end{equation}
and $\varphi_\mu$ its right-continuous inverse. Then the stopping time
\begin{equation}
  \label{eq:stoprefllevch3}
  \overline{T}_{\varphi_\mu}=\inf\big\{t>0:Y_t= \varphi_\mu(L_t)\big\}
\end{equation}
is a.s.\ finite and solves the Skorokhod embedding problem for $Y$, i.e.\ $Y_{\overline{T}_{\varphi_\mu}}\sim\mu$.
\end{proposition}
For measures without atoms the above proposition was derived independently by Pistorius \cite{pistorius_skoro}. The drawback of the above solution is that the stopping time $\overline{T}_{\varphi_\mu}$ does not need to be minimal.
Indeed, consider $\mu$ given as the law of $Y_R$ for $R=\inf\{t: Y_t\in \{1,2\}\}$. Then naturally $Y_R\sim Y_{\overline{T}_{\varphi_\mu}}$ and $R\le \overline{T}_{\varphi_\mu}$. Moreover, for large class of processes $X$, $\p(R<\overline{T}_{\varphi_\mu})>0$ which contradicts minimality of $\overline{T}_{\varphi_\mu}$. However minimality for discontinuous processes is  a hard issue and we do not have any feasible criterion, like in Brownian motion case (cf.~Cox \cite{cox_min}), to decide whether a stopping time is minimal. It might be thus rational to consider other ways of expressing the idea that a stopping time is small imposing, for example, some integrability conditions\footnote{Recall that a typical dummy solution to the Skorokhod embedding problem (cf.\ \cite[Sec.~2.3]{genealogia}) has infinite expectation.}. As shown by Pistorius \cite{pistorius_skoro} if $(X_t)$ drifts to $+\infty$ and $\mu$ is integrable, $\int_0^\infty s d\mu(s)<\infty$, then $\overline{T}_{\varphi_\mu}$ is integrable under the usual assumption: $D_\mu(\infty)=\int_0^\infty W(s)/W_+'(s) d\mu(s)<\infty$.

We turn now to the signed case and embedding for spectrally negative L\'evy process $(X_t)$. To develop an embedding for $(X_t)$ itself we need to redefine the sign of an excursion. We put $sgn(\epsilon)=\lim_{s\to 0}\epsilon(s)/|\epsilon(s)|$. As $X$ has only negative jumps, a negative excursion stays always negative, however an excursion which we call positive can either stay always positive or became negative from some moment on. We also redefine the extrema functional (\ref{eq:extremach3}) via $\overline{M}(\epsilon)(t)=sgn(\epsilon)\sup\{|\epsilon(s)|:s\le t\land V(\epsilon), sgn(\epsilon(s))=sgn(\epsilon)\}$, that is we only count positive extremum of a positive excursion. So equipped we can present a solution to the Skorokhod embedding for spectrally negative L\'evy processes.
\begin{proposition}
Let $(X_t)$ be a spectrally negative recurrent L\'evy process for which $\{0\}$ is regular and instantaneous and let $\mu$ be a non-atomic probability measure on $\re$ such that $\mu(\re_-)>0$ and $\mu(\re_+)>0$, and (\ref{eq:suppFch3}) and (\ref{eq:cond_regch3}) hold for $F=\overline{M}$. The the stopping time
\begin{equation}\label{eq:stopspecneglevych3}
    \overline{T}_{\varphi_-,\varphi_+}=\inf\big\{t>0: X_t=\varphi_+(L_t) \textrm{ or }\overline{M}_t \le X_t=-\varphi_-(L_t)\big\}
\end{equation}
where $\varphi_{+/-}$ are given by (\ref{eq:psi_expressch3})-(\ref{eq:varphisch3}), is a.s.\ finite and solves the Skorokhod embedding problem for $X$, i.e.\ $X_{\overline{T}_{\varphi_-,\varphi_+}}\sim \mu$.
\end{proposition}
The proposition follows from Theorem \ref{thm:mainch3} and its proof (presented in Section \ref{sec:solch3}) upon three simple observations: $L_{\overline{T}_{\varphi_-,\varphi_+}}=L_{T^{\overline{M}}_{\varphi_-,\varphi_+}}$, where $T^{\overline{M}}_{\varphi_-,\varphi_+}$ is defined in (\ref{eq:def_stop1ch3}), $sgn(X_{\overline{T}_{\varphi_-,\varphi_+}})=sgn(\overline{M}_{T^{\overline{M}}_{\varphi_-,\varphi_+}})$ and finally $X_{\overline{T}_{\varphi_-,\varphi_+}}\in \{-\varphi_-(L_{\overline{T}_{\varphi_-,\varphi_+}}),\varphi_+(L_{\overline{T}_{\varphi_-,\varphi_+}})\}$.

We choose to work with spectrally negative L\'evy processes, but it should be clear that we could equally well work with spectrally positive L\'evy processes.
The above solution demonstrates the flexibility of our approach and is of interest as there are very few explicit works in the discontinuous setup. It also complements well the solution of Bertoin and Le Jan \cite{MR93b:60173} which is rather adapted for symmetric L\'evy processes. However, we have to point out that the above solution has two disadvantages: the stopping time $\overline{T}_{\varphi_-,\varphi_+}$ is in general not minimal and the characteristic measure of $\overline{M}$ may be quite hard to compute.

\section{Atomic measures}
\label{sec:atomsch3}
We now turn our attention to probability measures $\mu$ with atoms. The reason for developing so far, in the signed case, only the setup of regular measures is simple: in the presence of atoms the formulae we obtain are much more involved. We can still describe an explicit solution to the Skorokhod embedding problem, however the functions $\psi_{+/-}$ will be given through an iteration procedure. We will not phrase our result as a theorem but rather develop it in steps. We treat only purely atomic measures.

Let $F$ be a functional as described in Section \ref{sec:funcch3} and $\mu$ a purely atomic probability measure, $\mu=\sum_{k=1}^n a_k\delta_{x_k}+\sum_{k=1}^m b_k\delta_{y_k}$, where $x_{k+1}<x_k<0$, $y_{k+1}>y_k>0$ and $a_k,b_k>0$, $\sum_{k=1}^n a_k+\sum_{k=1}^m b_k =1$, $n,m\le \infty$. Suppose that $n_F((-\infty,x_n])>0$ and $n_F([y_m,\infty))>0$. We will describe functions $\varphi_+$ and $\varphi_-$ such that the stopping time $T^F_{\varphi_-,\varphi_+}$ given by (\ref{eq:def_stopch3}) solves the Skorokhod embedding problem for $F$, i.e.\ $F_{T^F_{\varphi_-,\varphi_+}}\sim \mu$.

Naturally $\varphi_-:\re_+\to\{-x_1,\dots,-x_n\}$, $\varphi_+:\re_+\to\{y_1,\dots,y_m\}$ and $\varphi_{+/-}$ are piece-wise constant and non-decreasing. We thus have
\begin{equation}
  \label{eq:phiatomsch3}
  \varphi_-(z)=-x_k\quad \alpha_{k-1}\leq z< \alpha_k;\quad
  \varphi_+(y)=y_k\quad \beta_{k-1}\leq y< \beta_k,
\end{equation}
for some positive, finite, increasing sequences $(\alpha_k:1\leq k< n)$ and $(\beta_k:1\le k < m)$, and $\alpha_0=\beta_0=0$, $\alpha_n=\beta_m=\infty$.
Our aim is to determine the values $\alpha_k$, $\beta_k$ in terms of $n_F$ and $\mu$. For ease of notation we write $n_F(x)=n_F((-\infty,x])$ and $\overline{n}_F(y)=n_F([y,\infty))$ and $T=T^F_{\varphi_-,\varphi_+}$ the stopping time defined in the second line in (\ref{eq:def_stop1ch3}).

The formula in (\ref{eq:e_h(F)ch3}) is still valid and taking $h(z)=\mathbf{1}_{z=x_i}$ and $h(z)=\mathbf{1}_{z=y_j}$, we find
\begin{equation}
  \label{eq:genatomch3}
  \begin{split}
  a_i&=n_F(x_i)\int_{\alpha_{i-1}}^{\alpha_{i}}\p(L_T\ge l)dl,\quad 1\le i\le n \\
b_k&=\overline{n}_F(y_j)\int_{\beta_{j-1}}^{\beta_{j}}\p(L_T\ge l)dl,\quad 1\le j\le m
  \end{split}
\end{equation}
where $\p(L_T\ge l)$ is given through (\ref{eq:law_L_Tch3}). We will show how (\ref{eq:genatomch3}) identifies functions $\varphi_{+/-}$ uniquely. Put $\alpha_0=\beta_0=0$. Suppose that we know the values of $\alpha_0,\dots,\alpha_i$ and $\beta_0,\dots,\beta_j$, where $0\le i<n$, $0\le j< m$ are such that $\alpha_{i+1}\land \beta_{j+1}>\alpha_{i}\lor \beta_{j}$.
In particular we know the probabilities $\p(L_T\ge l)$ for $l\le \alpha_{i}\lor \beta_{j}$.
We will now show how to determine the smaller value of the two: $\alpha_{i+1}$, $\beta_{j+1}$. The procedure then repeats.

Assume for example that $\alpha_i\ge \beta_j$, the other case being similar. Then $\beta_{j+1}\leq \alpha_{i+1}$ if and only if
\begin{eqnarray}
\label{eq:condatomsch3}
&&  \int_{\alpha_i}^{\beta_{j+1}}\p (L_T\ge l)dl\le   \int_{\alpha_i}^{\alpha_{j+1}}\p (L_T\ge l)dl,\quad \textrm{if and only if}\nonumber\\
&&\frac{1}{\overline{n}_F(y_{j+1})}\Big[b_{j+1}-\overline{n}_F(y_{j+1})  \int_{\beta_j}^{\alpha_{i}}\p (L_T\ge l)dl\Big]\leq \frac{a_{i+1}}{n_F(x_{i+1})},\qquad
\end{eqnarray}
where we used (\ref{eq:genatomch3}) to obtain the last equivalence.
The last condition in (\ref{eq:condatomsch3}) is phrased uniquely in terms of known quantities. Thus we know how to determine which of the two values: $\alpha_{i+1}$, $\beta_{j+1}$ is smaller.
We note also that (\ref{eq:condatomsch3}) has a very natural interpretation which follows from basic properties of Poisson point processes. Namely, the value $\overline{n}_F(y_{j+1})  \int_{\beta_j}^{\alpha_{i}}\p (L_T\ge l)dl$ corresponds to the probability that we have stopped in $y_{j+1}$ for $L_T\le \alpha_i$. The value between parenthesis on the left hand side of (\ref{eq:condatomsch3}) corresponds simply to the weight that remains to be attributed to the atom in $y_{j+1}$. Therefore (\ref{eq:condatomsch3}) is a comparison between two ratios of the type: \emph{the weight attributed to the region} divided by \emph{the characteristic measure of the region}. For example, with $i=j=0$, we have $\beta_1\le \alpha_1$ if and only if $\frac{b_1}{\overline{n}_F(y_1)}\le \frac{a_1}{n_F(x_1)}$.

Suppose that $\beta_{j+1}\le \alpha_{j+1}$ or equivalently that (\ref{eq:condatomsch3}) is verified (the other case being symmetric). Then the value $\beta_{j+1}$ can be uniquely determined from (\ref{eq:genatomch3}). Namely, using (\ref{eq:law_L_Tch3}), we obtain the following equation for $\beta_{j+1}$
\begin{equation}
\label{eq:atomsstepcalcch3}
  \begin{split}
  b_{j+1}\,=\,&\overline{n}_F(y_{j+1})\bigg[\int_{\beta_i}^{\alpha_i}\p(L_T\ge l)dl
+\p(L_T\ge \alpha_i)\\&\cdot\int_{\alpha_{i}}^{\beta_{j+1}}dl \exp \Big\{-(l-\alpha_i)\big(n_F(x_{i+1})+\overline{n}_F(y_{j+1})\big)\Big\}\bigg] ,
  \end{split}
\end{equation}
which can be solved explicitly.

The conditions under which our procedure ends successfully determine the class of atomic probability measure that can be embedded using this method. Suppose that $n,m<\infty$. The construction works well until $i\leq n-1$ and $j\le m-1$. Suppose however that the procedure allowed us to define $\alpha_i$ for $0\le i\le n-1$ and $\beta_j$ for $0\le j\le k$ for some $k<m-1$. Then, as we need to have $\alpha_{n}=\infty$ the condition (\ref{eq:condatomsch3}) has to yield $\beta_j\le \alpha_n$ for $k\le j\le m-1$. Furthermore, for $k=m-1$ we need to have actually equality in the condition (\ref{eq:condatomsch3}).
These conditions extend in a natural way to the case when one of $n$ and $m$, or both, are infinite. Under such restrictions on the measure $\mu$ (which correspond to (\ref{eq:cond_regch3}) in the regular case), the procedure described in this section provides a solution to the Skorokhod embedding problem for $F$.\smallskip\\
\textbf{Example.} Let $n_F(dx)=dx/x^2$, $x\neq 0$ and $\mu=a\delta_{x_1}+b\delta_{y_1}+(1-a-b)\delta_{y_2}$ with $a,b>0$, $x_1<0<y_1<y_2$. With the notation used above $n=1$, $m=2$ so that $\alpha_1=\infty=\beta_2$. Thus, the condition (\ref{eq:condatomsch3}) has to yield $\beta_1<\alpha_1$ or equivalently $by_1<a|x_1|$, which is the first condition we impose on $\mu$. We then proceed to calculate $\beta_1$ using (\ref{eq:atomsstepcalcch3}):
\begin{equation}
  \label{eq:atomsexamplech3}
  \beta_1=\frac{|x_1|y_1}{y_1-x_1}\log\Big(\frac{x_1}{b(y_1-x_1)-x_1}\Big).
\end{equation}
To end the construction we have to require that (\ref{eq:condatomsch3}) yields $\alpha_1=\beta_2$ which is equivalent to $by_1+(1-a-b)y_2=ax_1$, that is to say $\mu$ is centered.
\smallskip\\
Following our methodology, one can also develop a solution to the Skorokhod embedding problem for $F$ for arbitrary measure $\mu$ with both regular and atomic components.
However, our description of such solution would be quite involved and we think there is no need to sketch it here, as it would bring little insight and the solution could hardly be qualified as \emph{explicit}.

\section{Proofs of the main results}
\label{sec:proofsch3}

In this last section we present the proofs of Theorems \ref{thm:mainch3} and \ref{thm:mainposch3}. We start with the former which is more technical and parts of which are then used in the proof of the latter.

\subsection{Proof of the embedding for signed functionals}
\label{sec:solch3}

In this section we prove Theorem \ref{thm:mainch3} and point out that the solution it presents is independent of the normalization of the local time and the It\^o measure of the underlying process $X$. Instead of just verifying that our embedding works we chose to present rather the complete reasoning which allows to obtain our solution. So, after proving that our stopping times are a.s.\ finite, we will ``pretend'' we do not know the explicit formulae for $\psi_+$ and $\psi_-$ and show how to discover them.

We start by calculating the law of $L_T$. For ease of notation, we denote the terminal value $F(\epsilon,V(\epsilon))$ simply by  $F(\epsilon)$. We have
{\setlength\arraycolsep{0.5pt}
\begin{eqnarray}
\label{eq:def_Nch3}
\p(L_{T^F_{\varphi_-,\varphi_+}}&&> l)=\p(T^F_{\varphi_-,\varphi_+}> \tau_l)\nonumber\\
&&=\p\Big(\textrm{on the time interval }[0,\tau_l]\textrm{ for every excursion }e_s,\; s\leq l,\nonumber\\
&&\quad\quad\; \textrm{the values of } F \textrm{ were between }-\varphi_-(s)\textrm{ and }\varphi_+(s)\Big)\nonumber\\
&&=\p\Big(\forall s\le l,\; F(e_s,V(e_s))\in\big(-\varphi_-(s),\varphi_+(s)\big)\Big)\nonumber\\
&&=\p\Big(\sum_{s\leq l}\mathbf{1}_{F(e_s)\notin(-\varphi_-(s),\varphi_+(s))}=0\Big)=\p(N_l=0),
\end{eqnarray}}
\noindent\ where the random variable $N_l=\sum_{s\leq l}\mathbf{1}_{F(e_s)\notin(-\varphi_-(s),\varphi_+(s))}$ is a Poisson variable with parameter
\begin{eqnarray}
  \label{eq:generaL_func_paramch3}
  &&\int_0^l n\Big(F(\epsilon)\notin(-\varphi_-(s),\varphi_+(s)\Big)ds\nonumber\\&=&\int_0^l n_F\Big((-\infty,-\varphi_-(s)]\cup[\varphi_+(s),+\infty)\Big)ds.
\end{eqnarray}
Thus, we obtain
\begin{equation}
  \label{eq:law_L_Tch3}
  \p(L_{T^F_{\varphi_-,\varphi_+}}> l)=\exp\bigg(-\int_0^l n_F\Big((-\infty,-\varphi_-(s)]\cup[\varphi_+(s),+\infty)\Big)ds\bigg).
\end{equation}
We note that the law of $L_{T^F_{\varphi_-,\varphi_+}}$ is absolutely continuous with respect to the Lebesgue measure. As $L_\infty=\infty$ a.s., the above gives us a convenient criterion for finiteness of our stopping time, namely $T^F_{\varphi_-,\varphi_+}<\infty$ a.s.\ if and only if the integral in (\ref{eq:generaL_func_paramch3}), with $l=\infty$, is infinite. We now prove the latter.

Recall that we assumed that $n_F((-\infty,x])>0$ and $n_F([y,\infty))>0$ for $a_\mu<x\leq 0\leq y< b_\mu$, and that (\ref{eq:cond_regch3}) holds. This ensures that the functions $\psi_+$ and $\psi_-$, given via (\ref{eq:psi_expressch3}) and (\ref{eq:psi2_expressch3}), are well defined. Let $\lambda=\psi_+(b_\mu)=\psi_+(\infty)$. We have $\psi_-(\infty)=\psi_-(-a_\mu)=\psi_+(D_\mu^{-1}(G_\mu(a_\mu)))=\psi_+(b_\mu)=\lambda$, where we used the assumption (\ref{eq:cond_regch3}) that $D_\mu(b_\mu)=G_\mu(a_\mu)$. We denote this last value by $c_\mu=D_\mu(b_\mu)$.

We need to calculate the integral in (\ref{eq:generaL_func_paramch3}) with $l=\infty$ and show that it is infinite. We have
\begin{eqnarray}
  \label{eq:firstch3}
  \int_0^\infty n_F([\varphi_+(s),\infty))ds&=&\int_0^{b_\mu}\frac{d\mu(s)}{1+\overline{\mu}(s)-\overline{\mu}(G^{-1}_\mu(D_\mu(s)))}\\
&=& \int_0^{c_\mu}\frac{n_F([D^{-1}_\mu (v),\infty))dv}{1+\overline{\mu}(D^{-1}_\mu (v))-\overline{\mu}(G^{-1}_\mu (v))},\nonumber
\end{eqnarray}
where the equalities follow with a change of variables from (\ref{eq:psi_expressch3}) and (\ref{eq:int_much3}). This is easy when $\mu$ has a positive density but is also true in the general setting. Indeed, since $\mu$ has no atoms, $\psi_+$ and $\psi_-$ are continuous and $\psi_{+/-}(\varphi_{+/-}(y))=y$. Jumps of $\varphi_{+/-}$ correspond to the level stretches of $\psi_{+/-}$, so that $d\psi_{+/-}$-- a.e.\ $\varphi_{+/-}(\psi_{+/-}(y))=y$. This justifies the first equality in (\ref{eq:firstch3}). For the second one, note that the functions $D_\mu$ and $G_\mu$ are constant only outside of the support of $\mu$, so that any $y>0$, $y\in supp(\mu)$, can be represented as $D^{-1}_\mu(u)$ and any $x<0$, $x\in supp(\mu)$, can be represented as $G^{-1}_\mu(v)$, for some $u,v$.
These remarks justify also the following derivation, based on (\ref{eq:psi2_expressch3}) and (\ref{eq:int_much3})
\begin{eqnarray}
  \label{eq:secondch3}
  \int_0^\infty n_F((-\infty,-\varphi_-(s)])ds&=&\int_{a_\mu}^0\frac{d\mu(s)}{1-\overline{\mu}(s)+\overline{\mu}(D^{-1}_\mu(G_\mu(s)))}\\
&=& \int_0^{c_\mu}\frac{n_F((-\infty,G^{-1}_\mu (v)])dv}{1+\overline{\mu}(D^{-1}_\mu (v))-\overline{\mu}(G^{-1}_\mu (v))}.\nonumber
\end{eqnarray}
Now observe that
\begin{eqnarray}
  \label{eq:diffgdch3}
d\Big(\overline{\mu}(D^{-1}_\mu (v))\Big)&=&-n_F([D^{-1}_\mu (v),\infty))dv\quad\textrm{and}\nonumber\\
d\Big(\overline{\mu}(G^{-1}_\mu (v))\Big)&=&n_F((-\infty,G^{-1}_\mu (v)])dv.
\end{eqnarray}
This allows us to calculate the desired integral in (\ref{eq:generaL_func_paramch3}). We have
\begin{eqnarray}
  \label{eq:finiteproofch3}
  &&\int_0^\infty n_F\Big((-\infty,-\varphi_-(s)]\cup[\varphi_+(s),+\infty)\Big)ds\quad\textrm{(using (\ref{eq:firstch3}) and (\ref{eq:secondch3}))}\nonumber\\
&&=\int_0^{c_\mu}\frac{n_F\Big((-\infty,G^{-1}_\mu (v)]\cup [D^{-1}_\mu (v),\infty)\Big)}{1+\overline{\mu}(D^{-1}_\mu (v))-\overline{\mu}(G^{-1}_\mu (v))}dv\quad\textrm{(using (\ref{eq:diffgdch3}))}\nonumber\\
&&=-\log\Big(1+\overline{\mu}(\infty)-\overline{\mu}(-\infty)\Big)=+\infty,
\end{eqnarray}
where we used the fact that $\overline{\mu}(D^{-1}_\mu(0))=\overline{\mu}(G^{-1}_\mu(0))$ and $D_\mu^{-1}(c_\mu)=\infty$, $G^{-1}_\mu (c_\mu)=-\infty$. We proved, by (\ref{eq:law_L_Tch3}), that $L_{T^F_{\varphi_-,\varphi_+}}<\infty$ a.s.\ and thus that $T^F_{\varphi_-,\varphi_+}<\infty$ a.s.
From (\ref{eq:finiteproofch3}) above, it can also be deduced that $\psi_+(b_\mu)+\psi_-(-a_\mu)=\infty$ and thus, as $\psi_+(b_\mu)=\psi_-(-a_\mu)$, we see that both are infinite.

We now turn to the proof of the embedding property announced in Theorem \ref{thm:mainch3}. Thanks to the property that the terminal value for $F_t$, for a given excursion, is either achieved on some interval or not achieved at all (see Section \ref{sec:funcch3}), we deduce instantly that $F_T\in\{-\varphi_-(L_T),\varphi_+(L_T)\}$.
From the properties of Poisson point processes, we see that conditionally on $\{L_T=l\}$, the respective probabilities that $F_T=-\varphi_-(L_T)$ or that $F_T=\varphi_+(L_T)$, are given by the proportions of the characteristic measures of appropriate regions, thus $\p(F_T=-\varphi_-(L_T)|L_T=l)=\frac{n_F((-\infty,-\varphi_-(l)])}{n_F((-\infty,-\varphi_-(l)]\cup[\varphi_+(l),+\infty))}$.\\ Let $h:\re\to\re_+$ be a bounded Borel function. We can then write
\begin{eqnarray}
\label{eq:e_h(F)ch3}
  \e h(F_T)&=&\e\bigg(\e\Big[h(F_T)\Big|L_T\Big]\bigg)
  =\e\bigg(\e\Big[h(-\varphi_-(L_T))\mathbf{1}_{F_T=-\varphi_-(L_T)}\nonumber\\ &&\hspace*{4cm}
+h(\varphi_+(L_T))\mathbf{1}_{F_T=\varphi_+(L_T)}\Big|L_T\Big]\bigg)\nonumber \\
  &=&\int_0^\infty \p(L_T\in dl)
\bigg[\frac{h(-\varphi_-(l))n_F((-\infty,\varphi_-(l)])}{n_F((-\infty,-\varphi_-(l)]\cup[\varphi_+(l),+\infty)\big)}\nonumber\\
&&\hspace*{2.8cm}+\frac{h(\varphi_+(l))n_F([\varphi_+(l),+\infty))}{n_F\big((-\infty,-\varphi_-(l)]\cup[\varphi_+(l),+\infty)\big)}\bigg]
\end{eqnarray}
The above formula, in Brownian setup and for the signed extrema functional (\ref{eq:extremach3}), was obtained by Jeulin and Yor \cite{MR83c:60110} (cf.\ Vallois \cite[Eq.~(2.3)]{MR86m:60200}).\\
On the other hand, since we want to have $F_T\sim \mu$, the above display (\ref{eq:e_h(F)ch3}) has to be equal to $\int_{\re} h(x)d\mu(x)$.
This has to be true for any bounded function $h$ and we will see that it will allow us to determine the functions $\varphi_-$ and $\varphi_+$.
Write $\psi_-$ and $\psi_+$ respectively for the inverses of $\varphi_-$ and $\varphi_+$, and assume that $\psi_+$, $\psi_-$ are continuous (recall that this is indeed our case since the measure $\mu$ in Theorem \ref{thm:mainch3} does not have any atoms).
Fix $y>0$ and put $h(z)=\mathbf{1}_{z\geq y}$. This yields
\begin{eqnarray}
\label{eq:ogon_pozch3}
 \overline{\mu}(y)&=&\int_{\psi_+(y)}^\infty\p(L_T\in dl)\bigg[\frac{n_F([\varphi_+(l),+\infty))}{n_F\big((-\infty,-\varphi_-(l)]\cup[\varphi_+(l),+\infty)\big)}\bigg]\nonumber\\
 &=&\int_{\psi_+(y)}^\infty dl\ n_F([\varphi_+(l),+\infty))\p(L_T\geq l),
\end{eqnarray}
where we differentiated (\ref{eq:law_L_Tch3}) to obtain $\p(L_T\in dl)$. Similarly, if we fix $x<0$ and put $h(z)=\mathbf{1}_{z\ge x}$, we obtain
\begin{eqnarray}
\label{eq:ogon_negch3}
  \overline{\mu}(x)&=&\p\big(L_T\leq \psi_-(-x)\big)
  + \int_{\psi_-(-x)}^\infty dl\ n_F([\varphi_+(l),+\infty))\p(L_T\geq l),
\end{eqnarray}
where we used the assumption that $\psi_-$ is continuous.
Assume that $f:\re_-\to \re_+$ and $g:\re_+\to \re_-$, given by
\begin{equation}
  \label{eq:def_f_gch3}
  f(x)=\varphi_+(\psi_-(-x))\quad\textrm{and}\quad g(y)=-\varphi_-(\psi_+(y))
\end{equation}
are well defined and finite for $a_\mu<x\leq 0\leq y<b_\mu$.\\
Recall the remarks between (\ref{eq:firstch3}) and (\ref{eq:secondch3}). In particular note that we can assume that  $f(g(y))=y$, $d\psi_+(y)$-- a.e., and $g(f(x))=x$, $d\psi_-(-x)$-- a.e.  \\ Let $y>0$ and differentiate (\ref{eq:ogon_pozch3}) to obtain
\begin{equation}
\label{eq:ogon_poz_diffch3}
    d\overline{\mu}(y)= -n_F([y,+\infty))\p(L_T\geq \psi_+(y))d\psi_+(y).
\end{equation}
Now it suffices to note that
\begin{eqnarray}
\label{eq:pomoc1ch3}
  \p\big(L_T<\psi_+(y)\big)&=&\p\big(L_T<\psi_-(-g(y))\big),\quad \textrm{by (\ref{eq:ogon_negch3})}\nonumber\\
&=&\overline{\mu}(g(y))-\int_{\psi_+(y)}^\infty dl\ n_F([\varphi_+(l),\infty))\p(L_T\geq l)\nonumber\\
&=&\overline{\mu}(g(y))-\overline{\mu}(y),\quad \textrm{using (\ref{eq:ogon_pozch3})}.
\end{eqnarray}
Combining (\ref{eq:ogon_poz_diffch3}) and (\ref{eq:pomoc1ch3}) above, we conclude that
\begin{equation}
  \label{eq:diff_psi_+ch3}
  d\psi_+(y)=\frac{-d\overline{\mu}(y)}{n_F\big([y,\infty)\big)\big(1+\overline{\mu}(y)-\overline{\mu}(g(y))\big)}.
\end{equation}
We will try to obtain similar expression starting with (\ref{eq:ogon_negch3}) instead of (\ref{eq:ogon_pozch3}). To this end fix $x<0$ and rewrite (\ref{eq:ogon_negch3}) using (\ref{eq:ogon_pozch3}) in the following way
\begin{equation}
\label{eq:mubar2ch3}
  \overline{\mu}(x)=1-\p(L_T\geq \psi_-(-x))+\overline{\mu}(f(x)),
\end{equation}
where we used the fact that $\psi_+(f(x))=\psi_-(-x)$. Note that we can assume that $\overline{\mu}(f(x))$ is continuous. Differentiating (\ref{eq:mubar2ch3}), through a reasoning similar to (\ref{eq:pomoc1ch3}), we obtain
\begin{eqnarray}
  d\overline{\mu}(x)&=&\Big(n_F((-\infty,x]) + n_F([f(x),+\infty))\Big)\nonumber\\&&\times\Big(1+\overline{\mu}(f(x))-\overline{\mu}(x)\Big)d\psi_-(-x)+d\overline{\mu}(f(x)),
\end{eqnarray}
where we used the fact that $d\psi_-$-- a.e.\ $\varphi_-(\psi_-(x))=-x$.
Taking $x=g(y)$, $y\geq 0$, in the above, yields
\begin{equation}
  \label{eq:diff_psi_+2ch3}
  d\psi_+(y)=\frac{d\overline{\mu}(g(y))-d\overline{\mu}(y)}{\Big(n_F((-\infty,g(y)]) + n_F([y,+\infty))\Big)\Big(1+\overline{\mu}(y)-\overline{\mu}(g(y))\Big)}.
\end{equation}
Comparing (\ref{eq:diff_psi_+ch3}) with (\ref{eq:diff_psi_+2ch3}) and simplifying the common terms we obtain finally
\begin{equation}
  \label{eq:diff_much3}
  \frac{d\overline{\mu}(y)}{n_F([y,+\infty))}=\frac{d\overline{\mu}(g(y))}{n_F((-\infty,g(y)])},\; y\geq 0.
\end{equation}
It is therefore natural to introduce the functions $D_\mu$ and $G_\mu$, defined in (\ref{eq:int_much3}), which we recall here
\begin{equation}
  \label{eq:int_mu2ch3}
  D_{\mu}(y)=\int_0^y \frac{d\mu(s)}{n_F([s,+\infty))}\quad\textrm{and}\quad G_{\mu}(x)=\int_{x}^0\frac{d\mu(s)}{n_F((-\infty,s])},
\end{equation}
for $y\ge 0$ and $x\le 0$. The functions $D_\mu$ and $G_\mu$ are continuous and increasing. Recall that their right-continuous inverses are denoted $D_\mu^{-1}$ and $G_\mu^{-1}$ respectively. The equality (\ref{eq:diff_much3}) reads $D_{\mu}(y)=G_{\mu}(g(y))$ or equivalently
\begin{equation}
  \label{eq:form_f_gch3}
  f(x)=D_{\mu}^{-1}(G_{\mu}(x)),\ x\leq 0,\quad\textrm{and}\quad g(y)=G_{\mu}^{-1}(D_{\mu}(y)),\ y\geq 0.
\end{equation}
Functions $f$ and $g$ were supposed to be well defined and we can now translate this assumption into conditions on $F$ and $\mu$. Namely, we need to have $D_\mu(y)<\infty$, $G_\mu(x)<\infty$ for $x<0<y$ and $D_\mu(\infty)=G_\mu(\infty)$.
The first condition means that for $x,y$ in the support of $\mu$ we need to have $n_F((-\infty,x])>0$ and $n_F([y,\infty))>0$, which we assumed in (\ref{eq:suppFch3}) and the second one is just (\ref{eq:cond_regch3}).

We are finally able to justify the explicit formulae for $\psi_+$ and $\psi_-$. Indeed, substituting the expression (\ref{eq:form_f_gch3}) for $g$ in (\ref{eq:diff_psi_+ch3}) and passing to the integral representation, we obtain (\ref{eq:psi_expressch3}), that is for $y\geq 0$,
\begin{equation}
  \label{eq:psi_express2ch3}
  \psi_+(y)=\int_0^y\frac{d\mu(s)}{n_F\Big([s,+\infty)\Big)\Big(1+\overline{\mu}(s)-\overline{\mu}\big(G^{-1}_{\mu}(D_\mu(s))\big)\Big)},
\end{equation}
and since $\psi_-(y)=\psi_+(f(-y))$, we have instantly
\begin{eqnarray}
  \label{eq:psi2_express2ch3}
  \psi_-(y)&=&\int_0^{D_\mu^{-1}(G_\mu(-y))}\frac{d\mu(s)}{n_F\Big([s,+\infty)\Big)\Big(1+\overline{\mu}(s)-\overline{\mu}\big(G^{-1}_{\mu}(D_\mu(s))\big)\Big)}\nonumber\\
&=&\int_{-y}^0\frac{d\mu(s)}{n_F\Big((-\infty,s]\Big)\Big(1+\overline{\mu}\big(D_\mu^{-1}(G_\mu(s))\big)-\overline{\mu}(s)\Big)},
\end{eqnarray}
where the second equality follows from (\ref{eq:diff_much3}) and (\ref{eq:form_f_gch3}).\\
It remains to prove that $T$ is minimal. Suppose that $S\leq T$ is a stopping time with $F_S\sim F_T$.
From the definition of $T$ in (\ref{eq:def_stop1ch3}) it follows that $0<F_S\leq \varphi_+(L_S)\leq \varphi_+(L_T)$ or $0\ge F_S\geq -\varphi_-(L_S)\ge -\varphi_-(L_T)$. We obtain for $\lambda,\gamma>0$
\begin{eqnarray}
  \label{eq:minimalitysign1ch3}
  \overline{\mu}(\lambda)&=&\e\mathbf{1}_{F_S\ge \lambda}\leq\e\Big[\mathbf{1}_{F_S\ge 0}\mathbf{1}_{\varphi_+(L_S)\ge\lambda}\Big]\leq\e\Big[\mathbf{1}_{F_S\ge 0}\mathbf{1}_{\varphi_+(L_T)\ge\lambda}\Big]
\nonumber\\&&= \overline{\mu}(\lambda)+\e\Big[\Big(\mathbf{1}_{F_S\ge 0}-\mathbf{1}_{F_T\ge 0}\Big)\mathbf{1}_{\varphi_+(L_T)\ge \lambda}\Big],\quad \textrm{and}\\
\mu((-\infty,-\gamma])&=&\e\mathbf{1}_{F_S\le \gamma}
\leq\e\Big[\mathbf{1}_{F_S\le 0}\mathbf{1}_{\varphi_-(L_S)\ge\gamma}\Big]\leq\e\Big[\mathbf{1}_{F_S\le 0}\mathbf{1}_{\varphi_-(L_T)\ge\gamma}\Big]\nonumber\\
&&=\mu((-\infty,-\gamma])+\e\Big[\Big(\mathbf{1}_{F_S\le 0}-\mathbf{1}_{F_T\le 0}\Big)\mathbf{1}_{\varphi_-(L_T)\ge \gamma}\Big],  \label{eq:minimalitysign2ch3}
\end{eqnarray}
from which it follows that $$\xi_\lambda:=\e\Big[\Big(\mathbf{1}_{F_S\ge 0}-\mathbf{1}_{F_T\ge 0}\Big)\mathbf{1}_{\varphi_+(L_T)\ge \lambda}\Big]\ge 0 \textrm{ and }
\eta_\gamma:=\e\Big[\Big(\mathbf{1}_{F_S\le 0}-\mathbf{1}_{F_T\le 0}\Big)\mathbf{1}_{\varphi_-(L_T)\ge \gamma}\Big]\ge 0.$$
Recall that $\mu$ has no atoms and therefore functions $\varphi_{+/-}$ are strictly increasing\footnote{This is not necessary for minimality of $T$ but simplifies the proof.}. It follows that $\xi_\lambda=-\eta_{\varphi_-(\psi_+(\lambda))}$ and thus $\xi_\lambda=\eta_\gamma=0$ for all $\lambda,\gamma>0$. This in turn signifies that we had equalities instead of inequalities everywhere in (\ref{eq:minimalitysign1ch3}) and (\ref{eq:minimalitysign2ch3}) and thus $L_T=L_S$, $F_S=\varphi_+(L_T)$ on $F_S>0$, and $F_S=-\varphi_-(L_T)$ on $F_S<0$. In consequence $T=S$.
This ends the proof of Theorem \ref{thm:mainch3}.

Recall that the excursion measure $n$ and the local time $L$ are uniquely determined up to a constant multiplicative factor. We argue now briefly that our results are independent of renormalization of $n$ and $L$. Suppose that we multiply the excursion measure by some constant $c$. The change $n\leadsto c\cdot n$ is directly translated into $n_F\leadsto c\cdot n_F$, which in turn gives $G_\mu\leadsto \frac{1}{c}G_\mu$ and $D_\mu\leadsto \frac{1}{c}D_\mu$. Such a transformation leaves unchanged both $D_\mu^{-1}(G_\mu(\cdot))$ and $G_\mu^{-1}(D_\mu(\cdot))$, thus leading to the change $\psi_+\leadsto \frac{1}{c}\psi_+$ and $\psi_-\leadsto \frac{1}{c}\psi_-$, which in turn translates into $\varphi_+(\cdot)\leadsto \varphi_+(c\cdot)$ and $\varphi_-(\cdot)\leadsto \varphi_-(c\cdot)$. Multiplying the excursion measure $n$ by a constant $c$ induces a multiplication of the local time at zero by the constant $\frac{1}{c}$, as can be readily seen from the relationship between the L\'evy measure of $\tau$ and the measure $n$ (cf.\ Bertoin \cite[Lemma IV.9]{MR1406564}), or from the fact that the measure $dn(\epsilon)dL_s$ is an invariant of excursion theory (which in turn follows easily from the compensation formula).
In consequence, the quantities $\varphi_+(L_t)$ and $\varphi_-(L_t)$ both stay unchanged, which implies that the stopping times given in (\ref{eq:def_stop1ch3}) stay unchanged as well. Thus our solution to the Skorokhod embedding problem for $F$ is independent of the normalization of the local time $L$ and the excursion measure $n$.

\subsection{Proof of the embedding for positive functionals}
\label{sec:proofposch3}
We now prove Theorem \ref{thm:mainposch3}. We start with its first part and consider $\mu$ with no atom in zero. We note that $\psi_\mu$ is increasing and $\psi_\mu(y)<\infty$ if $y\in supp(\mu)$. Indeed, we have then $\overline{\mu}(y)>0$ and $n_F([y,\infty))>0$ so that $\psi_\mu(y)\leq n_F([y,\infty))^{-1}(\frac{1-\overline{\mu}(y)}{\overline{\mu}}-\log (\overline{\mu}(y)))$. Secondly, note that $\psi_\mu(\infty)=\infty$. Indeed, for $0<\kappa<y<b_\mu$, we have $\psi_\mu(y+)\ge n_F([\kappa,\infty))^{-1}\log\big(\frac{\overline{\mu}(\kappa)}{\overline{\mu}(y+)}\big)\to\infty$ as $y\to b_\mu$.\\
We will use our previous study from Section \ref{sec:solch3}. Recall the notation from Section \ref{sec:mainch3}. If we put $\psi_-=0$, so that $\varphi_-=\infty$ then $T^F_{\varphi_\mu}=T^F_{\varphi_-,\varphi_+}$, displayed in (\ref{eq:def_stopch3}), with $\varphi_+=\varphi_\mu$. For simplicity we write $T$ for $T^F_{\varphi_\mu}$. The calculation of the law of $L_{T}$ given in (\ref{eq:law_L_Tch3}) is valid and with a change of variables it is seen that $L_T<\infty$ a.s.\ and thus the stopping time $T$ is a.s.\ finite for any probability measure $\mu$ with $\mu(\{0\})=0$. This follows also, as $\psi_\mu(\infty)=\infty$, from $\overline{\mu}(y)=\p(L_T\ge \psi_\mu(y))$ proved below.\\
Let $\nu$ be the non-atomic part of $\mu$, $d\nu(x)=\mathbf{1}_{\mu(\{x\})=0}d\mu(x)$ and $0<y_1<y_2<\dots$ the atoms of $\mu$.
We have to verify that for all $u\ge 0$, $\overline{\mu}(u)=\p(F_{T}\ge u)$. We have $\p(F_{T}\ge u)=\p(\varphi_\mu(L_T)\ge u)=\p(L_T\ge \psi_\mu(u))$ as $\psi_\mu$ and $\varphi_\mu$ are taken respectively left- and right- continuous. It suffices thus to show that $\overline{\mu}$ and $\p(L_T\ge \psi_\mu(u))$ have the same jumps and that they satisfy the same differential equation for all $u$ such that $\mu(\{u\})=0$ and  $u\in supp(\mu)$ (since the measures $d\psi_\mu$ and $d\mu$ have the same support). Observe that $u$ with $u\in supp(\mu)$, $\mu(\{u\})=0$ are precisely of the form $u=\varphi_\mu(y)$ with $\psi_\mu(\varphi_\mu(y))=y$.

Let $u=\varphi_\mu(y)$ such that $\mu(\{u\})=0$ so that $\psi_\mu(u)=y$. Note that $d\nu(u)=d\mu(u)$ and $dy=d(\psi_\mu(u))=d\mu(u)[\overline{\mu}(u)n_F([u,\infty))]^{-1}$. We have thus
\begin{equation}
  \label{eq:psiposregch3}
  \begin{split}
  d\Big(\p(L_T\ge \psi_\mu(u))\Big)&=
d\p(L_T\ge y)=-\p(L_T\ge y)n_F([\varphi_\mu(y),\infty))dy\\
&=\p(L_T\ge \psi_\mu(u))n_F([u,\infty))d\psi_\mu(u) \\
d\Big(\overline{\mu}(u)\Big)&=-\overline{\mu}(u)n_F([u,\infty))dy=-\overline{\mu}(u)n_F([u,\infty))d\psi_\mu(u),
  \end{split}
\end{equation}
which shows that the two functions $\overline{\mu}(\cdot)$ and $\p(L_T\ge \psi_\mu(\cdot))$ satisfy the same differential equation on the required set. It remains to show that $\mu(\{y_i\})=\p(\psi_\mu(y_i)\le L_T<\psi_\mu(y_i+))$, $i\ge 1$. Suppose we know this already for $i\le j-1$ for some $j\ge 2$. Combined with the discussion above, this yields $\overline{\mu}(y_j)=\p(L_T\ge \psi_\mu(y_j))$. Note also that for $y\in [\psi_\mu(y_j),\psi_\mu(y_j+))$, $\varphi_\mu(y)=y_j$ and let $\Delta \psi_\mu(y_j)=\psi_\mu(y_j+)-\psi_\mu(y_j)=\ln\big(\frac{\overline{\mu}(y_i)}{\overline{\mu}(y_i+)}\big)n_F([y_i,\infty))^{-1}$. By (\ref{eq:law_L_Tch3}) we have
\begin{eqnarray}
\label{eq:derivatomposch3}
 && \p(\psi_\mu(y_j)\le L_T<\psi_\mu(y_j+))=\p(L_T\ge \psi_\mu(y_j))-\p(L_T\ge \psi_\mu(y_j+)\nonumber\nonumber\\
&&=\p(L_T\ge \psi_\mu(y_j))\Big(1-\exp\big\{-\Delta\psi_\mu(y_j)n_F([y_j,\infty))\big\}\Big)\nonumber\nonumber\\
&&=\overline{\mu}(y_j)\Big(1-\frac{\overline{\mu}(y_j+)}{\overline{\mu}(y_j)}\Big)=\overline{\mu}(y_j)-\overline{\mu}(y_j+)=\mu(\{y_j\}),
\end{eqnarray}
which ends the proof of the embedding. It remains to verify that $T$ is minimal.
From the definition of $T$ in (\ref{eq:stopposch3}), as $\varphi_\mu$ is non-decreasing, it is clear that $\sup_{t\le T^F_{\varphi_\mu}}F_t=F_{T^F_{\varphi_\mu}}$. Let $S$ be a stopping time such that $F_S\sim F_T$ and $S\le T$. We have $\sup_{t\le S}F_t\le \sup_{t\le T^F_{\varphi_\mu}}F_t=F_{T^F_{\varphi_\mu}}$ and thus $0\le F_S\le F_T$ which, together with $F_T\sim F_S$, implies that $F_S=F_T$ a.s.\ and thus $S=T$ a.s. This ends the proof of the first part of Theorem \ref{thm:mainposch3}.

To finish the proof, we now argue the embedding for a probability measure $\mu$ with an atom in zero. It is an easy consequence of what we have obtained so far. Recall that we defined a new probability measure on $(0,\infty]$ transferring the atom in zero to infinity, $\tilde{\mu}=\mu-\varsigma (\delta_{0}-\delta_{\infty})$, where $\varsigma=\mu(\{0\})$. We thus have $\varsigma=\overline{\tilde{\mu}}(\infty)>0$ which implies $\psi_{\tilde{\mu}}(\infty)=\delta<\infty$. In turn, $\varphi_{\tilde{\mu}}(x)=\infty$ for $x\ge \delta$ and $\p(L_{T^F_{\varphi_{\tilde{\mu}}}}=\infty)=\p(L_{T^F_{\varphi_{\tilde{\mu}}}}\ge \delta)=\varsigma$, and $F_{T^F_{\varphi_{\tilde{\mu}}}}\mathbf{1}_{L_{T^F_{\varphi_{\tilde{\mu}}}}<\delta}\sim \tilde{\mu}|_{(0,\infty)}\sim \mu|_{(0,\infty)}$. It suffices now to observe that for $R=\inf\{t: L_t=\delta\}$ we have $F_R=0$ as the support of $dL_t$ is contained in the set of zeros of $(F_t)$, to conclude that $F_{R\land T^F_{\varphi_{\tilde{\mu}}}}\sim \mu$.

\section{Closing remarks}
\label{sec:endch3}
When starting the present work, we set ourselves two goals. The first goal was to develop, following the methodology outlined by Ob\l\'oj and Yor \cite{obloj_yor}, a \emph{two-sided} (and thus non-randomized) solution to the Skorokhod embedding problem for Az\'ema's martingale. The second goal was to understand the essence of such a solution, as well as of the solution presented in \cite{obloj_yor} and generalize them to the most abstract setting in which they still work. Put differently, we wanted to understand what were the necessary properties to be imposed on considered processes, and what was the convenient abstract framework for phrasing the solutions. We hope, we have achieved our goals with the presentation of general solution to the Skorokhod embedding for class of processes, functionals of Markovian excursions, given in Section \ref{sec:mainch3}. We stress that our solution is given in terms of the law of the functional under It\^o's measure.
These laws have been studied now for over 30 years.
In a lot of important cases they are well known, as Sections \ref{sec:azemach3} and \ref{sec:pozch3} witness. We also dispose of various general methods to investigate them (cf.\ the link with Krein's string theory indicated in Section \ref{sec:funcch3}). This makes, we hope, our solution workable and useful.

We close this paper with a final remark about the functionals of excursions we have used throughout the paper. We saw in Section \ref{sec:disconch3} that we can, in certain cases, develop explicit embeddings using discontinuous functionals. In fact, we could also generalize the class of admissible functionals introduced in Section \ref{sec:funcch3} in a different direction.
We restricted the setup of this paper to functionals $F$ which satisfy $|F(\epsilon,V(\epsilon))|>0$ for any generic non-trivial excursion $\epsilon$. This was equivalent with requiring that the sets of zeros of $(F_t)$ and $(X_t)$ coincide, which in turn implied that our random times $T^F$, displayed in (\ref{eq:def_stop1ch3}), were stopping times in the natural filtration of $F$. It should be clear however, that
if one agrees to work with stopping times relative to a larger filtration, as the natural filtration of $X$, then our study generalizes instantly to functionals $F$ which can be equivalent to zero at some excursions.
An important example of such functional is given by the sojourn time of a one-dimensional diffusion $X$, during an excursion, above a certain level $\lambda$: $F^\lambda(\epsilon,t)=\int_0^{t\land V(\epsilon)}\mathbf{1}_{\epsilon_s\ge \lambda}ds$. The process $F^\lambda_t=\int_{g_t}^t\mathbf{1}_{X_s\ge \lambda}ds$ is closely linked with the sojourn time for $X$, $\int_0^t\mathbf{1}_{X_s\ge \lambda}ds$. The latter is an additive functional studied by number of authors (cf.\ K\"uchler \cite{MR857226,MR901531}, Truman and Williams \cite{MR1182602}) and the characteristic measure $n_{F^\lambda}$ can be calculated.
\bigskip\\
{\large \textbf{Acknowledgement}}\smallskip\\ Author wants to express his gratitude towards Marc Yor for his constant help and support.\newpage
\bibliographystyle{plain}
\bibliography{../../bibliografia}
\end{document}